\documentclass[10pt]{article}
\usepackage{mypackages}
\usepackage{mycommands}

\title{Optimal explicit stabilized postprocessed $\tau$-leap method for the simulation of chemical kinetics}

\author{Assyr Abdulle\thanks{\email{assyr.abdulle@epfl.ch}} }
\author{Lia Gander\thanks{\email{lia.gander@usi.ch}, currently at the Euler Institute, Universit\`a della Svizzera italiana, Lugano, Switzerland.} }
\author{Giacomo Rosilho de Souza\thanks{\email{giacomo.rosilhodesouza@epfl.ch}}}
\affil{ANMC, Institute of Mathematics, \'Ecole Polytechnique F\'ed\'erale de Lausanne, Lausanne, Switzerland}

\begin{document}

\maketitle

\begin{abstract}
The simulation of chemical kinetics involving multiple scales constitutes a modeling challenge (from ordinary differential equations to Markov chain) and a computational challenge (multiple scales, large dynamical systems, time step restrictions). In this paper we propose a new discrete stochastic simulation algorithm: the postprocessed second kind stabilized orthogonal $\tau$-leap Runge--Kutta method (PSK-$\tau$-ROCK). In the context of chemical kinetics this method can be seen as a stabilization of Gillespie's explicit $\tau$-leap combined with a postprocessor.
The stabilized procedure allows to simulate problems with multiple scales (stiff), while the postprocessing procedure allows to approximate the invariant measure (e.g. mean and variance) of ergodic stochastic dynamical systems. We prove stability and accuracy of the PSK-$\tau$-ROCK. 
Numerical experiments illustrate the high reliability and efficiency of the scheme when compared to other $\tau$-leap methods.
\end{abstract}


\vspace{0.5cm}

\section{Introduction}\label{sec:intro}

The modeling of kinetic chemical processes involves multiple chemical species with different population's size and reacting time-scales. 
A typical ordinary differential equation (ODE) model for the simulation of such problems are the reaction rate equations (RRE), but this model is accurate only at the thermodynamic limit (i.e. when the populations size and system volume tend to infinity, but the concentrations remain constant). In contrast, for systems with small populations, as living cells, discrete and stochastic modeling is required to capture the correct kinetics. 
Assuming proper mixing and thermal equilibrium a discrete dynamical system in the form of a Markov process as well as its accompanying master equation, the chemical master equation (CME), can be derived for the evolution of the probability density function of a chemical system \cite{Gil92,Qua67}.
The stochastic simulation algorithm (SSA) \cite{Gil76,Gil77} gives an exact method to compute samples from the distribution of the CME. However, while very easy to implement, the SSA can become overwhelmingly slow due to the presence of multiple scales in the system (stiffness) and presence of large and small populations size, leading to reactions that fire extremely often.

By fixing a step size (or leap time) $\tau$ encompassing several reactions, Gillespie proposes the $\tau$-leap method \cite{Gil01}. This approximate procedure lumps together the reactions that would occur in a time lapse $\tau$ and fires them simultaneously. If the number of reactions fired in a time $\tau$ is large then the $\tau$-leap scheme approximates the Euler--Maruyama method for the chemical Langevin equation (CLE) and in the thermodynamic limit this latter scheme approaches the explicit Euler method for the RRE \cite{Gil01,Hig08}. 

In this paper we focus on $\tau$-leap methods. As the reactions fire at disparate time-scales and the systems are typically stiff, the standard explicit $\tau$-leap method \cite{Gil01} faces stability issues \cite{CPR04}. Furthermore, even when the stability conditions are met, amplification properties due to explicitness of the scheme prevent to capture the correct statistics of the process. Implicit $\tau$-leap schemes \cite{CGP07,CPR04,HMT17,RPC03} in contrast usually do not have stability constraints. However, for ergodic dynamical systems, these schemes generally fail in capturing the exact statistics of the system. Hence, unless the fast processes are resolved, both explicit and implicit $\tau$-leap methods fail to correctly integrate ergodic dynamical systems (see \cite{LAE08} for a similar discussion in the context of stochastic differential equations (SDEs) driven by diffusion processes). We mention further the trapezoidal $\tau$-leap method, that is accurate in sampling the invariant measure for linear equations, but might fail for nonlinear problems \cite{LAE08}.
Very recently, in \cite{RKV19}, a split step scheme generalizing the $\theta$-method is introduced. The method is accurate in sampling the invariant measure of the process, however at each step it requires the solution of two nonlinear systems and an optimization problem for the scheme's coefficients.

Apart from the aforementioned implicit or explicit Runge--Kutta-like methods, several hybrid schemes making use of different models and levels of granularity exist in the literature. Such multirate (or multiscale) methods exploit the multiscale nature of chemical reaction systems, which often consist of multiple reactions firing at disparate time-scales. Roughly speaking, most of such schemes divide the reactions (or species) into fast and slow ones. Then, the fast dynamics are resolved by making use of a quasiequilibrium assumption and the slow terms are integrated employing larger step sizes --- see \cite{CGP05a,ELV05b,HVH19,HaR02,HZC97,HAL12,RaA03}. In this paper we do not assume that the system is clearly separable into fast and slow dynamics, therefore multirate methods are not discussed in what follows.

We now briefly describe explicit stabilized methods, that are the basis building blocks of our new scheme. In the ODE context, explicit stabilized methods are a compromise between explicit and implicit schemes. No linear algebra solutions are needed as for implicit methods, while quadratic growth (with the number of stages) of their stability domains allows for much better stability properties than classical explicit methods.
Well known explicit stabilized methods are the Runge--Kutta--Chebyshev (RKC) methods \cite{SSV98,HoS80,VHS90}, the DUMKA methods \cite{Leb94,LeM94,Med98} and the Runge--Kutta orthogonal Chebyshev (ROCK) \cite{Abd02,AbM01} methods (note that the first-order RKC and ROCK scheme coincide). More recently, the first-order RKC (or ROCK) scheme has been extended to SDEs, yielding the S-ROCK family \cite{AbC07,AbC08,AbL08,Blu15} and higher order extensions in \cite{AVZ13b}. 
For mean-square stable problems the S-ROCK scheme introduced in \cite{AbL08} represents an important improvement over the Euler--Maruyama method thanks to its improved stability properties (it does however not preserve the optimal stability domain of the first-order RKC method). However, for non mean-square stable problems and for problems with too large variance the efficiency of the S-ROCK scheme deteriorates. 
Starting from the same S-ROCK scheme as in \cite{AbL08}, the authors in \cite{AHL10} derive the $\tau$-ROCK method for equations driven by discrete noise; however, this method inherits the same issues as the S-ROCK method of \cite{AbL08}.
The SK-ROCK scheme \cite{AAV18} is an improvement over the previous S-ROCK method for SDEs, this scheme has an optimal stability domain's size growing quadratically with the number of stages and is second-order accurate in sampling the invariant measure of a class of ergodic SDEs. 

The main contribution of this paper is the design and analysis of a new $\tau$-leap method for stiff chemical systems. 
Inspired by the SK-ROCK scheme we propose here the PSK-$\tau$-ROCK leap method for chemical kinetics. This method has several desirable properties:
\begin{itemize}%[label=\roman*)]
	\item it is fully explicit avoiding any linear algebra computations and is as easy to implement as the explicit $\tau$-leap method;
	\item it has an extended and optimal stability domain growing quadratically with the number of function evaluations, avoiding any step size restriction as the standard explicit $\tau$-leap method;
	\item thanks to a postprocessing technique adapted from \cite{AAV18,Vil15} it shows remarkable properties in sampling correct statistics of non mean-square stable chemical systems, even when fast reactions are not resolved. 
\end{itemize}
We analyze the accuracy and stability properties of the scheme and its long-time dynamics for ergodic dynamical systems.
The efficiency and accuracy of the new scheme are illustrated through a sequence of numerical experiments, where we also compare the method against other $\tau$-leap schemes as the implicit $\tau$-leap and the trapezoidal $\tau$-leap method \cite{CPR04,RPC03}. 

The rest of the paper is organized as follows. In \cref{sec:ssa_tau} we give an introduction to the SSA and $\tau$-leap method, in \cref{sec:psktrock} we introduce the PSK-$\tau$-ROCK method and provide a detailed pseudocode. The accuracy and stability analysis of the scheme is given in \cref{sec:analysis}, while in \cref{sec:numexp} we provide the numerical examples. Conclusions are found in \cref{sec:conclu}.

\section{The SSA and the \texorpdfstring{$\tau$}{tau}-leap method}\label{sec:ssa_tau}
In this section we briefly recall the modeling of a well stirred chemical reaction system at thermal equilibrium and introduce the SSA and $\tau$-leap scheme.

\paragraph{A model for chemical reaction systems.}
Consider a chemical system composed by $N$ species (of molecules) $S_1,\ldots,S_N$ which interact in $M$ reactions, denoted $R_1,\ldots,R_M$. We are interested in the number of molecules of each specie in an instant of time $t$. We denote by $\bX(t)=(X_1(t),\ldots, X_N(t))^\top$ the state vector, where $X_j(t)\in \Nb$ is the number of molecules of specie $S_j$ at time $t$. 
Each reaction $R_j$ is characterized by the propensity function $a_j(\bx)$ and the state-change vector $\bnu_j$. Given a state $\bx\in \Nb^N$ and an infinitesimal time $\dif t$, the quantity $a_j(\bx)\dif t$ is the probability that reaction $R_j$ fires within $\dif t$ units of time. The state-change vector $\bnu_j$ describes the change in state $\bx$ when reaction $R_j$ is fired, i.e. reaction $R_j$ has the effect of changing the state vector from $\bx$ to $\bx+\bnu_j$. We will denote by $a(\bx)=(a_1(\bx),\ldots,a_M(\bx))^\top$ the vector of propensity functions and by $\bnu=(\bnu_1,\ldots,\bnu_M)$ the stoichiometric matrix.

\begin{example}\label{ex:MM}
	We provide here an illustrative example of the above description. To do so, we consider the famous Michaelis--Menten system describing the mechanism of enzymatic catalysis. The model consists in four species: a substrate $S_1$, an enzyme $S_2$, a complex enzyme-substrate $S_3$ and the product $S_4$. The three reactions may be written as
	\begin{equation}\label{eq:MM}
	S_1+S_2 \stackrel{c_1}{\rightharpoonup} S_3, \qquad\qquad
	S_3 \stackrel{c_2}{\rightharpoonup} S_1+S_2,\qquad\qquad
	S_3 \stackrel{c_3}{\rightharpoonup} S_4+S_2.
	\end{equation}
	The state vector is $\bX(t)=(X_1(t),X_2(t),X_3(t),X_4(t))^\top$ and represents the number of molecules of each specie $S_1,S_2,S_3,S_4$ at time $t$. If the first reaction fires, the value of $X_3(t)$ is increased by one molecule and $X_1(t),X_2(t)$ are decreased by one molecule each, hence the state vector is updated as $\bX(t)+\bm{\nu}_1$, where $\bm{\nu}_1=(-1,-1,1,0)^\top$. In the same manner we define $\bm{\nu}_2=(1,1,-1,0)^\top$ and $\bm{\nu}_3=(0,1,-1,1)^\top$. The propensity function $a_1(\bx)$ is the probability that the first reaction takes place within one unit of time and is given by $a_1(\bx)=c_1 x_1 x_2$, where the product $x_1x_2$ is the number of possible distinct combinations of $S_1,S_2$ molecules and $c_1$ is the probability that given two reactants $S_1,S_2$ the reaction actually fires. Similarly, $a_2(\bx)= c_2x_3$ and $a_3(\bx)=c_3x_3$. 
	
	There are also other types of reaction which are not listed in \cref{eq:MM}, for instance $2S_n\stackrel{c_j}{\rightarrow} S_m$ and $3S_n\stackrel{c_j}{\rightarrow} S_m$, whose propensity functions are $a_j(\bx)=c_jx_n(x_n-1)/2!$ and $a_j(\bx)=c_jx_n(x_n-1)(x_n-2)/3!$, respectively, and their structure follows from a combinatoric argument.
\end{example}

Given the vector of propensity functions $a(\bx)$ and the stoichiometric matrix $\bnu$ the system evolves following two simple rules \cite{Gil76}.
\begin{enumerate}[label=\roman*)]
	\item Given a state vector $\bX(t)$ at time $t$, in an infinitesimal time $\dif t$ the reactions $R_j$ are independent and the probability that reaction $R_j$ fires is given by $a_j(\bX(t))\dif t$. \label{rules:first}
	\item If $R_j$ fires the system is updated as $\bX(t+\dif t) = \bX(t)+\bm{\nu}_j$.
\end{enumerate}

\paragraph{The Stochastic Simulation Algorithm.}
From \ref{rules:first}, Gillespie \cite{Gil76} derived a probability density function from which we can sample a random pair $(\tau,j)$, where $\tau$ is the waiting time until the next reaction and $j$ is the index of the next reaction. This is the core of the stochastic simulation algorithm (SSA), given by:

\begin{enumerate}[label=\arabic*)]
\item Sample the waiting time $\tau$ from an exponentially distributed random variable with rate $a_0(\bX(t))=\sum_{j=1}^M a_j(\bX(t))$.\label{ssa:first}
\item Sample the next reaction $j$ from an $M$ point random variable, where index $j$ has probability $a_j(\bX(t))/a_0(\bX(t))$.
\item Update the state vector as $\bX(t+\tau)=\bX(t)+\bnu_j$ and the time as $t\leftarrow t+\tau$.
\item Return to \ref{ssa:first}, unless a stopping criteria is satisfied.
\end{enumerate}
The most important property of the SSA is that it is exact in sampling the statistics of the system. However, if there is at least one reaction with high probability of firing, then $a_0(\bX(t))$ will be large and the waiting time $\tau$ will likely be very short. Therefore, the SSA will use an excessively large number of time steps and become practically unreasonably expensive.

\paragraph{The \texorpdfstring{$\tau$}{tau}-leap method.}
The $\tau$-leap method \cite{Gil01} speeds the simulation by fixing a step size $\tau$ and firing all the reactions that occur within time $\tau$ simultaneously. This leaping strategy leads to a good approximation of the SSA if the so-called \emph{leap condition} is satisfied: the propensity functions $a_j(\bx)$ must not change appreciably in the time interval $[t,t+\tau]$.

First, suppose that in the time interval $[t,t+\tau]$ the propensity functions $a_j(\bx)$ are constant and thus the reaction events are independent. Under this assumption, the number of times that reaction $R_j$ fires in the time interval $[t,t+\tau]$ is described by a Poisson random variable with rate $a_j(\bx)\tau$, that we denote as $\mathcal{P}_j(a_j(\bx)\tau)$. Hence, under the leap condition, the $\tau$-leap scheme
\begin{equation}\label{eq:tauleap}
	\bX_{n+1} = \bX_n + \sum_{j=1}^M \bm{\nu}_j \mathcal{P}_j(a_j(\bX_n)\tau)
\end{equation}
is a good approximation to the SSA, where $\bX_n$ is an approximation of $\bX(\tn)$ with $\tn=n\tau$. We note that in order to satisfy the leap condition the reactants population cannot be too small, otherwise a few reactions change considerably the number of reactants and thus the propensity functions change substantially as well.

Since the mean (and variance) of $\mathcal{P}_j(a_j(\bx)\tau)$ is $a_j(\bx)\tau$, it is useful to decompose the right-hand side of \cref{eq:tauleap} in a drift term and a zero-mean noise term:
\begin{equation}\label{eq:tauleapEM}
	\bX_{n+1} = \bX_n +\tau f(\bX_n) + Q(\bX_n,\tau),
\end{equation}
where
\begin{equation}\label{eq:deffQtauleap}
	f(\bx) =  \sum_{j=1}^M \bm{\nu}_j a_j(\bx), \qquad\qquad Q(\bx,\tau) =  \sum_{j=1}^M \bm{\nu}_j (\mathcal{P}_j(a_j(\bx)\tau)-a_j(\bx)\tau).
\end{equation}
We note that \cref{eq:tauleapEM} is very similar to the Euler--Maruyama scheme for SDEs, where the diffusion is replaced by the zero-mean discrete noise $Q(\bx,\tau)$.

For stiff chemical systems, the approximation \cref{eq:tauleapEM} can face severe step size $\tau$ restrictions to be stable \cite{CPR04}. Using implicit time-stepping can cure stability issues at the expense of solving nonlinear problems. But implicit methods might fail to capture the correct statistics of a chemical system, in form of mean and variance, due to damping introduced by implicitness.

\section{The PSK-\texorpdfstring{$\tau$}{tau}-ROCK method}\label{sec:psktrock}
In this section we introduce the PSK-$\tau$-ROCK scheme. This explicit stabilized $\tau$-leap method is composed of:
\begin{enumerate}[label=\roman*)]
	\item a time-marching scheme (denoted SK-$\tau$-ROCK) for the computation of approximate solutions $\bX_n$;
	\item a postprocessing procedure (denoted P) used to improve the accuracy of $\bX_n$ whenever needed, usually only at the very last time step.
\end{enumerate}
In \cref{sec:algo} we define the time-marching scheme SK-$\tau$-ROCK while in \cref{sec:postalgo} we motivate and introduce the postprocessing procedure P. The combination of the SK-$\tau$-ROCK time-marching scheme with the postprocessor P yields the PSK-$\tau$-ROCK scheme. In \cref{sec:pseudocode} we provide a detailed pseudocode for the PSK-$\tau$-ROCK scheme and discuss some implementation details.

Considering a test problem, we will show in \cref{sec:analysis} that the PSK-$\tau$-ROCK scheme has an optimal stability domain growing quadratically with the number of stages and thanks to the postprocessing procedure accurate sampling of the process' statistics is achieved. The properties shown on the test problem in \cref{sec:analysis} are verified numerically in \cref{sec:numexp} on more involved problems.

\subsection{The SK-\texorpdfstring{$\tau$}{tau}-ROCK step}\label{sec:algo}
Let $\tau$ be the step size, $\varepsilon\geq 0$ be the damping parameter and $\beta=2-4\varepsilon/3$; typically $\varepsilon=0.05$. We denote by $\rho$ the spectral radius of the Jacobian of $f$ evaluated in $\bX_n$, with $f$ as in \cref{eq:deffQtauleap}, and let the number of stages $s\in\Nb$ satisfy $\tau\rho\leq \beta s^2$. The SK-$\tau$-ROCK step, of size $\tau$, is given by 
\begin{align}\label{eq:sktaurock}
\begin{split}
\bK_0&=\bX_n,\\
\bK_1 &= \bK_0+\mu_1\tau f(\bK_0+\nu_1 Q(\bK_0,\tau))+\kappa_1 Q(\bK_0,\tau),\\
\bK_j&= \nu_j \bK_{j-1}+\kappa_j \bK_{j-2}+\mu_j\tau f(\bK_{j-1}), \quad j=2,\ldots,s,\\
\bX_{n+1}&=\bK_s,
\end{split}
\end{align}
where $f,Q$ are given in \cref{eq:deffQtauleap}. The coefficients $\mu_j,\nu_j,\kappa_j$, for $j=1,\ldots,s$, are as follows. We let
\begin{equation}\label{eq:defozou}
\oz=1+\varepsilon/s^2,\qquad\qquad \ou=T_s(\oz)/T_s'(\oz),
\end{equation}
where $T_s(x)$ is the Chebyshev polynomial of the first kind of degree $s$, defined by
\begin{equation}
	T_0(x)=1,\qquad\qquad T_1(x)=x,\qquad\qquad T_j(x)=2xT_{j-1}(x)-T_{j-2}(x) \quad j\geq 2.
\end{equation}
Finally, we define $\mu_1=\ou/\oz$, $\nu_1=s\ou/(2\oz)$\footnote{Our definition of $\nu_1$ is slightly different than in \cite{AAV18}. The motivation of this modification will be given in \cref{sec:analysis}.}, $\kappa_1=s\ou/\oz$ and, for $j=2,\ldots,s$,
\begin{equation}\label{eq:defmnk}
	\mu_j=2\ou T_{j-1}(\oz)/T_j(\oz),\qquad \nu_j=2\oz T_{j-1}(\oz)/T_j(\oz),\qquad \kappa_j=-T_{j-2}(\oz)/T_j(\oz).
\end{equation}

In \cref{eq:sktaurock}, only one evaluation of the drift term $f$ is required for accuracy, while the additional $s-1$ evaluations are used to increase stability. Indeed, as we will see in \cref{sec:analysis}, the SK-$\tau$-ROCK step involves the first and second kind shifted Chebyshev polynomials, that are instrumental to obtain optimal stability domains. 
The parameter $\varepsilon$ in \cref{eq:defozou} is called damping parameter. For $\varepsilon=0$ the stability domain of the method \cref{eq:sktaurock} will have a finite number of points $z_i$ along the negative real axis for which the absolute value of the stability function is exactly one. This is avoided setting $\varepsilon>0$. Also, introduction of damping is essential to study the ergodic properties of the numerical scheme (see \cref{sec:analysis} for details).

The main difference with respect to the previous $\tau$-ROCK scheme \cite{AHL10} is that here the noise term is put at the beginning of the iteration and therefore it is stabilized by the drift. In the reversed $\tau$-ROCK scheme, also introduced in \cite{AHL10}, the noise is as well put at the beginning of the iteration but with different parameters $\nu_1=1$ and $\kappa_1=0$, yielding in an overly damped noise.

\subsection{The postprocessing procedure}\label{sec:postalgo}
In chemical reactions, one is often interested in the stationary state of a given system. Hence, an algorithm capable of capturing the invariant measure of the system is of considerable interest. Therefore, we propose here a postprocessing procedure for the SK-$\tau$-ROCK time-marching scheme introduced in \cref{sec:algo}, which allows to considerably improve its accuracy when applied to non mean-square stable problems. We stress that the postprocessor is applied only when higher accuracy is required and it is not needed to advance the solution in time. However, before introducing the postprocessing procedure for chemical kinetics we briefly motivate it recalling the postprocessors' theory for linear SDEs. 

\paragraph{Postprocessors for linear SDEs.}
Postprocessors are since long employed to increase the accuracy of numerical solutions to ODEs \cite{But69}. However, a postprocessors framework for ergodic SDEs has been only recently proposed in \cite{Vil15}. We recall here the ideas developed in \cite{Vil15} but restricted to the very particular case of the Ornstein--Uhlenbeck process. Consider the SDE
\begin{equation}\label{eq:OU}
	\dif X(t)=\lambda X(t)\dif t+\sigma \dif W(t),\qquad\qquad X(0)=X_0,
\end{equation}
where $X_0\in\Rb$ is deterministic, $X(t)\in\Rb$, $W(t)$ is a Wiener process and $\lambda,\sigma\in\Rb$ with $\lambda<0$. The exact solution $X(t)$ is Gaussian with 
\begin{equation}
	\lim_{t\to\infty} \exp(X(t))=0,\qquad \lim_{t\to\infty}\var(X(t))=\frac{\sigma^2}{2|\lambda|}.
\end{equation}
Applying a Runge--Kutta method to \cref{eq:OU} yields $X_{n+1}=A(z)X_n+B(z)\sqrt{\tau}\sigma\xi_n$, where $z=\tau\lambda$ and $\xi_n\sim\mathcal{N}(0,1)$. Using recursion we deduce that, if $|A(z)|<1$,
\begin{equation}
	\lim_{n\to\infty} \exp(X_n)=0,\qquad \lim_{n\to\infty}\var(X_n)=\frac{\sigma^2}{2|\lambda|}R(z), \quad\mbox{with}\quad R(z)=\frac{-2zB(z)^2}{1-A(z)^2}.
\end{equation}
Therefore, the numerical method has order $r_1$ for the invariant measure (i.e. $|\lim_{n\to\infty}\var(X_n)-\var(X_\infty)|=\bigo{z^{r_1}}
$) if, and only if, $R(z)=1+\bigo{z^{r_1}}$ as $z\to 0$. However, higher order is easily achieved applying a postprocessing procedure. Indeed, applying the postprocessor $\overline X_n = C(z)X_n+D(z)\sqrt{\tau}\sigma\xi_n$ yields
\begin{equation}
	\lim_{n\to\infty} \exp(\overline X_n)=0,\qquad \lim_{n\to\infty}\var(\overline X_n)=\frac{\sigma^2}{2|\lambda|}(C(z)^2R(z)-2zD(z)^2)
\end{equation}
and therefore higher order $r_2>r_1$ is achieved choosing $C(z),D(z)$ such that $C(z)^2R(z)-2zD(z)^2=1+\bigo{z^{r_2}}$ as $z\to 0$.

\paragraph{The postprocessing procedure.}
Based on the ideas developed in \cite{Vil15} for the Ornstein--Uhlenbeck process, we define here the postprocessor for the SK-$\tau$-ROCK time-marching scheme \cref{eq:sktaurock}. Up to our knowledge, in the literature of chemical kinetics no such postprocessors have been used.

In order to obtain higher accuracy for the invariant measure of the system at a certain time $\tn=n\tau$, the postprocessor
\begin{equation}\label{eq:postproc}
\overline{\bX}_{n} =  \bX_{n} +\alpha\, Q(\bX_{n},\tau),
\end{equation}
with
\begin{equation}\label{eq:defalpha}
\alpha = \frac{1}{2}\sqrt{\frac{\ou}{\oz}},
\end{equation}
is employed. 
We stress that the PSK-$\tau$-ROCK scheme does not need to compute \cref{eq:postproc} at each time step but only whenever higher accuracy is required.

Due to the damping properties of the SK-$\tau$-ROCK steps (see \cref{thm:sktaurock} below), the variance of the numerical solution $\bX_n$ in \cref{eq:sktaurock} is smaller than the exact variance. Adding the random variable $\alpha Q(\bX_n,\tau)$ in \cref{eq:postproc} allows to increase the variance of the numerical solution, yielding in a better approximation.

\subsection{The PSK-\texorpdfstring{$\tau$}{tau}-ROCK method: the algorithm}\label{sec:pseudocode}
The PSK-$\tau$-ROCK method, thus, advances the solution in time using the SK-$\tau$-ROCK scheme \cref{eq:sktaurock,eq:defozou,eq:defmnk} and applies the cheap postprocessing step \cref{eq:postproc,eq:defalpha}, whenever higher accuracy for the invariant measure is needed. We summarize in this section such method by providing a detailed pseudocode in \cref{algo:pseudocode}.

The input parameters of \cref{algo:pseudocode} are the initial value $\bX_0$, the step size $\tau$, the end time $T$ and the drift and compensated Poisson noise terms $f(\bx)$ and $Q(\bx,\tau)$, respectively, which are defined in \cref{eq:deffQtauleap}. The output is the postprocessed numerical solution $\overline{\bX}_N$, which is an approximation to the exact solution $\bX(T)$, with $T=N\tau$. The procedure for computing the method's coefficients at \cref{line:callcoeff} of \cref{algo:pseudocode} is given in Function \ref{algo:coeff}($s,\varepsilon$) below.
We conclude the section with a few comments on \cref{algo:pseudocode}.
\begin{itemize}
	\item Approximation of the spectral radius at \cref{line:rho} is very cheap if performed with nonlinear power methods \cite{Lin72,Ver80}. In our experience those methods usually converge with at most two function evaluations, see \cref{tab:MMtime,tab:nlrevreac,tab:genloopcpu}. It is good practice to store the eigenvector associated to the largest eigenvalue and use it as starting guess for the next call to the nonlinear power method.
	\item We emphasize that \cref{algo:pseudocode} has low memory requirement as it needs three stage vectors $\bK_{-1,0,1}$ only, disregarding the size of $s$. Moreover, \cref{line:memswap} has zero cost if performed by simply swapping memory addresses.
	\item It is common to replace $s$ by $s+1$ after \cref{line:defs}, this enlarges the stability domain and ensures stability of the method even if the spectral radius $\rho$ increases within one time step.
	\item The call to \texttt{Coefficients}($s,\varepsilon$) at \cref{line:callcoeff} is needed only if the number of stages $s$ changes from one time step to the next. This does not happen too frequently. 
	\item \Cref{algo:pseudocode} and Function \ref{algo:coeff}($s,\varepsilon$) can be merged. Indeed the computation of coefficients $\mu_j,\nu_j,\kappa_j$ for $j=2,\ldots,s$ can be done inside the for loop beginning at \cref{line:for} of \cref{algo:pseudocode}, avoiding the execution of a for loop exclusively for the coefficients' definition. However, this does not improve significantly the performance unless the number of stages changes frequently.
	\item Finally, \cref{line:ou} in Function \ref{algo:coeff}($s,\varepsilon$) has negligible cost if the values of $\ou$ are precomputed and stored in table. 
\end{itemize}

\IfStandalone
{\begin{algorithm}[H]}{\begin{algorithm}[ht]}
		\DontPrintSemicolon
		\LinesNumbered
		\SetKwInOut{Input}{Input}
		\SetKwInOut{Output}{Output}
		\SetKwFunction{Coefficients}{Coefficients}
		\Input{$\bX_0$, $\tau$, $T$, $f(\bx)$, $Q(\bx,\tau)$}
		\Output{$\overline{\bX}_N$}
		%		\BlankLine
		$\varepsilon=0.05$ and $\beta =2-4\varepsilon/3$\;
		$t_0=0$\;
		\While{$t_0<T$}
		{
			Approximate the spectral radius $\rho$ of the Jacobian of $f$ evaluated at $\bX_0$\;\nllabel{line:rho}
			Let $s\in\Nb$ be the smallest integer satisfying $\tau\rho\leq\beta s^2$\; \nllabel{line:defs}
			$(\mu,\nu,\kappa,\alpha)=\Coefficients(s,\varepsilon)$\;\nllabel{line:callcoeff}
			$\bK_0=\bX_0$\;
			$\bK_1=\bK_0+\mu_1\tau f(\bK_0+\nu_1 Q(\bK_0,\tau))+\kappa_1 Q(\bK_0,\tau)$\;
			\For{$j\leftarrow 2$ \KwTo $s$}{ \nllabel{line:for}
				$\bK_{-1}=\bK_0$ and $\bK_0=\bK_1$\; \nllabel{line:memswap}
				$\bK_1= \nu_j \bK_{0}+\kappa_j \bK_{-1}+\mu_j\tau f(\bK_{0})$\;
			}
			$\bX_{0} = \bK_1$\;
			$t_0\leftarrow t_0+\tau$\;
		}
		$\overline{\bX}_{N}=\bX_{0}+\alpha\, Q(\bX_{0},\tau)$\;
		\caption{The PSK-$\tau$-ROCK method}
		\label{algo:pseudocode}
	\end{algorithm}
	
	\IfStandalone
	{\begin{function}[H]}{\begin{function}[ht]}
			\DontPrintSemicolon
			\LinesNumbered
			\SetKwInOut{Input}{Input}
			\SetKwInOut{Output}{Output}
			\Input{$s$ and $\varepsilon$}
			\Output{Method's coefficients $\mu_j,\nu_j,\kappa_j$ for $j=1,\ldots,s$}
			$\oz=1+\varepsilon/s^2$,
			$\ou=T_s(\oz)/T_s'(\oz)$\;\nllabel{line:ou}
			$\alpha = \sqrt{\ou/\oz}/2$\;
			$\mu_1=\ou/\oz$,
			$\nu_1=s\ou/(2\oz)$,
			$\kappa_1=s\ou/\oz$\;
			$T_0(\oz)=1$, $T_1(\oz)=\oz$\;
			\For{$j\leftarrow 2$ \KwTo $s$}{
				$T_j(\oz) = 2\oz T_{j-1}(\oz)-T_{j-2}(\oz)$\;
				$\mu_j=2\ou T_{j-1}(\oz)/T_j(\oz)$, $\nu_j=2\oz T_{j-1}(\oz)/T_j(\oz)$, $\kappa_j=-T_{j-2}(\oz)/T_j(\oz)$\;
			}
			\caption{Coefficients($s$,$\,\varepsilon$)}
			\label{algo:coeff}
		\end{function}

\IfStandalone
{
	\bibliographystyle{siam}
	\bibliography{../../../../../../LaTeX/library}
}{}

\section{Accuracy and stability analysis}\label{sec:analysis}
We analyze here the long-time accuracy and stability of the PSK-$\tau$-ROCK scheme on a model problem: the reversible isomerization reaction
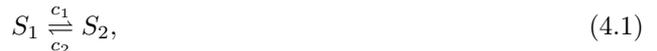
\begin{equation}\label{eq:reviso}
S_1  \underset{c_2}{\stackrel{c_1}{\rightleftharpoons}} S_2,
\end{equation}
which was first introduced in \cite{CPR04} and plays the role of test equation. This is a reversible system and therefore the number of molecules is constant, i.e. $X_1(t)+X_2(t)=X^T$. As a consequence specie $S_2$ can be neglected, we consider only specie $S_1$ and denote $X(t)=X_1(t)$. The system is described by the two reactions
\begin{equation}\label{eq:reacreviso}
a_1(x) = c_1 x, \qquad a_2(x) = c_2(X^T-x), \qquad \bm{\nu}_1 = -1, \qquad \bm{\nu}_2 = 1.
\end{equation}
Note that, for this particular model, $\lambda=-(c_1+c_2)$ is the only eigenvalue of the Jacobian of $f$, with $f$ as in \cref{eq:deffQtauleap}. Hence, the spectral radius of the Jacobian of $f$ is $\rho=|\lambda|$. In what follows, $\lambda$ represents the stiffness of the system. 

Problem \cref{eq:reviso} has a stationary state $X_\infty$ with a binomial distribution $B(n,p)$, where $n=X^T$ and $p=c_2/|\lambda|$ \cite{Gil02}. Hence,
\begin{equation}\label{eq:invmes}
\exp(X_\infty) = \frac{c_2}{|\lambda|}X^T, \qquad\qquad \var(X_\infty)=\frac{c_1c_2}{|\lambda|^2}X^T.
\end{equation}
Note that if $c_2=0$, i.e. \cref{eq:reviso} is not reversible, then $\exp(X_\infty) =\var(X_\infty)=0$ and the problem is considered to be mean-square stable. 

\begin{definition}
	A numerical method applied to \cref{eq:reviso} is said to have absolutely stable mean and variance if, and only if, $\exp(X_n)$ and $\var(X_n)$ remain bounded as $n\to\infty$. 	
	Moreover, if $c_2=0$ and $\lim_{n\to\infty}\exp(X_n)=0$ and $\lim_{n\to\infty}\var(X_n)=0$ then the method is said to be mean-square stable.
\end{definition}

\paragraph{Preliminary results.}
In the foregoing analysis we will need the polynomials
\begin{equation}
A_s(z)=\frac{T_s(\oz+\ou z)}{T_s(\oz)}, \qquad\qquad B_s(z)=\frac{U_{s-1}(\oz+\ou z)}{U_{s-1}(\oz)}\left( 1+\frac{\ou}{2\oz }z\right),
\end{equation}
where $U_n(x)$ is the Chebyshev polynomial of the second kind of degree $n$ defined recursively by
\begin{equation}
	U_0(x)=1,\qquad\qquad U_1(x)=2x,\qquad\qquad U_j(x)=2xU_{j-1}(x)-U_{j-2}(x) \quad j\geq 2
\end{equation}
and $\oz,\ou$ are defined in \cref{eq:defozou}; therefore, $A_s,B_s$ depend implicitly on $\varepsilon$. In the next lemma we collect known results about the stability polynomial $A_s(z)$ \cite{HoS80}. In particular, observe as the stability domain of $A_s(z)$ grows quadratically with $s$.
\begin{lemma}\label{lemma:resultsAs}
Let $s\in\Nb$, $\varepsilon\geq 0$ and $\ell_s^\varepsilon=2\oz/\ou$. Then $\ell_s^\varepsilon\geq \beta s^2$, with $\beta=2-4\varepsilon/3$, and the stability polynomial $A_s(z)$ satisfies $|A_s(z)|\leq 1$ for all $z\in [-\ell_s^\varepsilon,0]$. 
If $\varepsilon>0$ then $\ell_s^\varepsilon>\beta s^2$ and $|A_s(z)|< 1$ for all $z\in (-\ell_s^\varepsilon,0)$. 
	%The same results hold for $B_s(z)$.
\end{lemma}

In what follows, we denote by SK-$\tau$-ROCK the $\tau$-leap scheme defined by \cref{eq:sktaurock}, but where no postprocessor is applied, i.e. it is the PSK-$\tau$-ROCK scheme without postprocessing \cref{eq:postproc}.
We state now a theorem showing the long-term stability and accuracy properties of this SK-$\tau$-ROCK method \cref{eq:sktaurock} applied to \cref{eq:reviso}, it is the starting point to study the accuracy of the PSK-$\tau$-ROCK scheme. We recall that for \cref{eq:reviso} $\rho=|\lambda|$.
\begin{theorem}\label{thm:sktaurock}
	If $\tau\rho\leq \beta s^2$ and $\varepsilon>0$ the mean and variance of the SK-$\tau$-ROCK scheme \cref{eq:sktaurock} applied to \cref{eq:reviso} are absolutely stable and satisfy 
	\begin{equation}\label{eq:sktautockstat}
	\lim_{n\to\infty} \exp(X_n) = \exp(X_\infty), \qquad\qquad \lim_{n\to\infty} \var(X_n)  =c_{s}(z) \var(X_\infty),
	\end{equation}
	with
	\begin{equation}\label{eq:defc}
	 c_s(z) = \frac{-2 z B_s(z)^2}{1-A_s(z)^2},
	\end{equation}
	$z=\tau\lambda$ and $ \exp(X_\infty)$, $\var(X_\infty)$ as in \cref{eq:invmes}.
\end{theorem}
\begin{proof}
	From \cref{eq:reacreviso,eq:deffQtauleap}, we deduce that for the test equation \cref{eq:reviso} it holds
	\begin{equation}\label{eq:fQtest}
		f(x)= -c_1 x+c_2(X^T-x)=\lambda x+c_2 X^T, \qquad\qquad Q(x,\tau) = N_2(x,\tau)-N_1(x,\tau),
	\end{equation}
	where 
	\begin{equation}
		N_1(x,\tau)= \mathcal{P}_1(c_1 x\tau)-c_1x\tau, \qquad\qquad N_2(x,\tau)=\mathcal{P}_2(c_2(X^T-x)\tau)-c_2(X^T-x)\tau.
	\end{equation}
	%For the stages $K_j$ in \cref{eq:sktaurock} we consider the change of variables $\tK_j=K_j+c_2X^T/\lambda$ and observe that $f(K_j)=\lambda \tK_j$. 
	Replacing $f$ as in \cref{eq:fQtest} into \cref{eq:sktaurock} yields
	\begin{equation}
		K_1 = K_0+ \mu_1\tau(\lambda K_0+c_2X^T)+(\mu_1\nu_1\tau\lambda+\kappa_1)Q(K_0,\tau)
	\end{equation}
	and thus, considering the change of variables $\tK_j=K_j+c_2X^T/\lambda$,
	\begin{equation}\label{eq:pertsktaurocka}
		\tK_1 = \tK_0+ \mu_1 z\tK_0+(\mu_1\nu_1 z+\kappa_1)Q(K_0,\tau),
	\end{equation}
	where $z=\tau\lambda$. Using $\nu_j+\kappa_j=1$ for $j=2,\ldots,s$ we also obtain
	\begin{equation}\label{eq:pertsktaurockb}
		\tK_j = \nu_j\tK_{j-1}+\kappa_j \tK_{j-2}+\mu_j z \tK_{j-1},\quad j=2,\ldots,s.
	\end{equation}
	It is shown recursively \cite{AAV18,VHS90} that	
	\begin{equation}\label{eq:stabsktaurocktmp2}
		\tK_{s} = A_s(z)\tK_0 + B_s(z)Q(X_n,\tau).
	\end{equation}
	Let $\tX_n =X_n+c_2 X^T/\lambda$, 
	%as $X_{n+1}=K_s$ and $X_n=K_0$ 
	thus $\tX_{n+1}=\tK_s$ and $\tX_n=\tK_0$. 
	Since $\exp(Q(X_n,\tau))=0$ then, from \cref{eq:stabsktaurocktmp2},
	\begin{equation}\label{eq:exptXn}
		\exp(\tX_{n+1}) = A_s(z)\exp(\tX_n)= A_s(z)^{n+1}\exp(\tX_0).
	\end{equation}
	Using \cref{lemma:resultsAs}, $\varepsilon>0$ and $|z|\leq \beta s^2$ it follows $|A_s(z)|< 1$, implying
	\begin{equation}\label{eq:limitXn}
		\lim_{n\to\infty}\exp(X_{n}) =\lim_{n\to\infty}\exp(\tX_{n}) -\frac{c_2}{\lambda}X^T =  \frac{c_2}{|\lambda|}X^T = \exp(X_\infty).
	\end{equation}
	For the variance we use the law of total variation $\var(X)=\exp(\var(X|Y))+\var(\exp(X|Y))$, where $X,Y$ are two random variables, and we obtain
	\begin{align}\label{eq:stabsktaurocktmp3}
		\begin{split}
			\var(Q(X_n,\tau))&= c_1\exp(X_n)\tau + c_2(X^T-\exp(X_n))\tau
			=(c_1-c_2)\exp(\tX_n)\tau-2\frac{c_1c_2}{\lambda}X^T\tau \\
			&= (c_1-c_2)A_s(z)^n\exp(\tX_0)\tau-2\frac{c_1c_2}{\lambda}X^T\tau .
		\end{split}
	\end{align}
	Therefore, from $\var(X_n)=\var(\tX_n)$, the law of total variation and \cref{eq:stabsktaurocktmp2}, follows
%	\begin{align}\label{eq:varXnexact}
%		\begin{split}
%			\var(X_{n+1}) &= A_s(z)^2\var(X_n)+B_s(z)^2\var(Q(X_n,\tau))\\
%			&= A_s(z)^{2(n+1)}\var(X_0) + B_s(z)^2\sum_{k=0}^n A_s(z)^{2k}\var(Q(X_{n-k},\tau))\\
%			&= A_s(z)^{2(n+1)}\var(X_0) + (c_1-c_2)\exp(\tX_0)\tau B_s(z)^2\sum_{k=0}^n A_s(z)^{n+k}\\
%			&\quad -2\frac{c_1c_2}{\lambda}X^T\tau B_s(z)^2\sum_{k=0}^n A_s(z)^{2k}\\
%			&= A_s(z)^{2(n+1)}\var(X_0) + (c_1-c_2)\exp(\tX_0)\tau B_s(z)^2A_s(z)^{n}\frac{1-A_s(z)^{n+1}}{1-A_s(z)}\\
%			&\quad -2\frac{c_1c_2}{\lambda}X^T\tau B_s(z)^2\frac{1-A_s(z)^{2(n+1)}}{1-A_s(z)^2}
%		\end{split}
%	\end{align}
\begin{equation}
		\var(X_{n+1}) = A_s(z)^2\var(X_n)+B_s(z)^2\var(Q(X_n,\tau)).
\end{equation}
From \cref{eq:stabsktaurocktmp3} we deduce, by recursion,
\begin{align}\label{eq:varXnexact}
	\var(X_{n+1})&= A_s(z)^{2(n+1)}\var(X_0) + (c_1-c_2)\exp(\tX_0)\tau B_s(z)^2A_s(z)^{n}\frac{1-A_s(z)^{n+1}}{1-A_s(z)}\\
	&\quad -2\frac{c_1c_2}{\lambda}X^T\tau B_s(z)^2\frac{1-A_s(z)^{2(n+1)}}{1-A_s(z)^2}
\end{align}
	and thus
	\begin{equation}
		\lim_{n\to\infty}\var(X_{n}) =  -2\frac{c_1c_2}{\lambda}X^T\tau \frac{B_s(z)^2}{1-A_s(z)^2} =  \frac{-2\tau\lambda B_s(z)^2}{1-A_s(z)^2}\frac{c_1c_2}{\lambda^2}X^T= \frac{-2 z B_s(z)^2}{1-A_s(z)^2}\var(X_\infty). \tag*{\qedhere}
	\end{equation}
\end{proof}

\begin{corollary}
	If $\tau\rho\leq \beta s^2$ and $\varepsilon>0$, the SK-$\tau$-ROCK scheme \cref{eq:sktaurock} applied to \cref{eq:reviso} with $c_2=0$ is mean-square stable.
\end{corollary}

We recall that $A_s(z),B_s(z)$ implicitly depend on $\varepsilon$ (through $\oz,\ou$), therefore $c_s(z)$ in \cref{eq:defc} depends on $\varepsilon$ too. The condition on the damping parameter, $\varepsilon>0$, is needed to ensure that $|A_s(z)|<1$; indeed, for $\varepsilon=0$ there is a finite number of points $z_i$ along the negative real axis where $|A_s(z_i)|=1$, which might trigger instabilities and also preclude the computation of the limit in \cref{eq:limitXn}.

We see in \cref{thm:sktaurock} that the SK-$\tau$-ROCK scheme (and thus also the PSK-$\tau$-ROCK) has the same stability condition as the RKC and SK-ROCK schemes, i.e. $\tau\rho \leq \beta s^2$, with $\beta=2-4\varepsilon/3$, $\varepsilon>0$.
Hence, up to damping, the SK-$\tau$-ROCK method has an optimal stability domain growing quadratically with the number of function evaluations, recall that the optimal stability domain grows as $2s^2$.
In contrast, for the $\tau$-ROCK and Rev-$\tau$-ROCK method \cite{AHL10}, the damping parameter increases with $s$ and therefore the size $\ell_s^\varepsilon$ of the stability domain decreases with $s$. We compare in \cref{fig:comp_eps} the damping parameters of the SK-$\tau$-ROCK and the $\tau$-ROCK methods. In \cref{fig:comp_ell} we compare the size $\ell_s^\varepsilon$ of the stability domains, note as $\ell_s^\varepsilon$ grows faster for SK-$\tau$-ROCK. We do not illustrate the results for the Rev-$\tau$-ROCK scheme as the damping parameter depends also on the spectral radius of the Jacobian of $f$, see \cite[Section 4.3]{AHL10}.
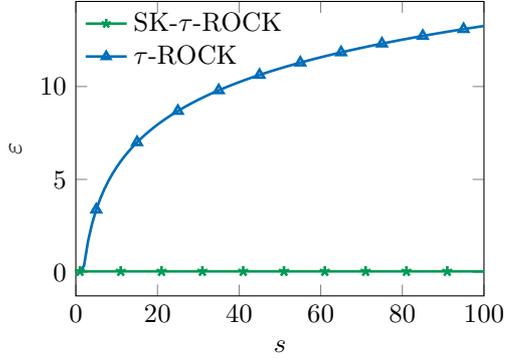
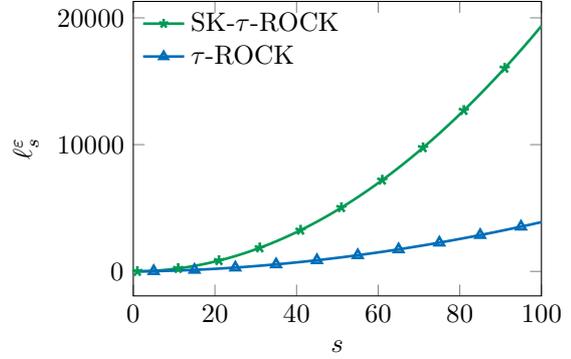
\begin{figure}
	\begin{center}
		\begin{subfigure}[t]{\subfigsize\textwidth}
			\centering
				\begin{tikzpicture}[every pin/.style={fill=white}]
					\begin{axis}[height=\plotheightd\textheight,width=\plotwidthd\textwidth,xmax=100,xmin=0,legend columns=1,legend style={draw=\legendboxdraw,fill=\legendboxfill,at={(0,1)},anchor=north west},legend cell align={left},
						xlabel={$s$}, ylabel={$\varepsilon$},label style={font=\normalsize},tick label style={font=\normalsize},legend image post style={scale=\legendmarkscale},legend style={nodes={scale=\legendfontscale, transform shape}},grid=none,scaled ticks=false]
						\addplot[color=colorone,line width=\plotlinewidth pt,mark=\markone,mark repeat=10,mark phase=1] table [x=slist,y=sktreps,col sep=comma] 
						{data/text/stab_dom_comp.csv};\addlegendentry{SK-$\tau$-ROCK};
						\addplot[color=colortwo,line width=\plotlinewidth pt,mark=\marktwo,mark repeat=10,mark phase=5] table [x=slist,y=treps,col sep=comma] 
						{data/text/stab_dom_comp.csv};\addlegendentry{$\tau$-ROCK};
					\end{axis}
				\end{tikzpicture}
			\caption{Displaying the damping parameter $\varepsilon$, with respect to $s$, of the SK-$\tau$-ROCK and $\tau$-ROCK methods.}
			\label{fig:comp_eps}
		\end{subfigure}\hfill
		\begin{subfigure}[t]{\subfigsize\textwidth}
			\centering
			\begin{tikzpicture}[every pin/.style={fill=white}]
				\begin{axis}[height=\plotheightd\textheight,width=\plotwidthd\textwidth,xmax=100,xmin=0,legend columns=1,legend style={draw=\legendboxdraw,fill=\legendboxfill,at={(0,1)},anchor=north west},legend cell align={left},
					xlabel={$s$}, ylabel={$\ell_s^\varepsilon$},label style={font=\normalsize},ticklabel style={font=\normalsize}, y tick label style={/pgf/number format/std,/pgf/number format/set thousands separator={}},legend image post style={scale=\legendmarkscale},legend style={nodes={scale=\legendfontscale, transform shape}},grid=none,scaled ticks=false]
					\addplot[color=colorone,line width=\plotlinewidth pt,mark=\markone,mark repeat=10,mark phase=1] table [x=slist,y=sktrell,col sep=comma] 
					{data/text/stab_dom_comp.csv};\addlegendentry{SK-$\tau$-ROCK};
					\addplot[color=colortwo,line width=\plotlinewidth pt,mark=\marktwo,mark repeat=10,mark phase=5] table [x=slist,y=trell,col sep=comma] 
					{data/text/stab_dom_comp.csv};\addlegendentry{$\tau$-ROCK};
				\end{axis}
			\end{tikzpicture}
			\caption{Displaying the stability domain size $\ell_s^\varepsilon$, with respect to $s$, of the SK-$\tau$-ROCK and $\tau$-ROCK methods.}
			\label{fig:comp_ell}
		\end{subfigure}
	\end{center}
	\caption{Comparison between the damping parameter $\varepsilon$ and the stability domain size $\ell_s^\varepsilon$ of the SK-$\tau$-ROCK and the 	$\tau$-ROCK method. For SK-$\tau$-ROCK $\varepsilon=0.05$, for $\tau$-ROCK $\varepsilon$ depends on $s$.}
	\label{fig:stabdom}
\end{figure}
For problems where the fast variables have non trivial invariant measure is was noted \cite[Example 3]{AHL10} that for the $\tau$-ROCK method values of $\varepsilon$ as high as $700$ were needed. However, as it is shown in \cite[Lemma 4.4]{AbL08}, for large damping values the stability domain grows only linearly with function evaluations, yielding in no gain with respect to the standard $\tau$-leap method. Moreover, a too large damping entails an excessively large number of stages, which introduce non negligible round-off errors \cite{VHS90}.

From \cref{eq:sktautockstat} we deduce that the SK-$\tau$-ROCK scheme is able to capture the exact expectation; in contrast, we observe in \cref{fig:ceps} that in general $c_s(z)<1$ and therefore the variance is under estimated. In the next section we show that the postprocessor of \cref{sec:postalgo} increases the accuracy for the variance and therefore the PSK-$\tau$-ROCK method is more accurate than the not postprocessed SK-$\tau$-ROCK method.
\begin{figure}
	\begin{center}
		\begin{subfigure}[t]{\subfigsize\textwidth}
			\centering
			\begin{tikzpicture}[every pin/.style={fill=white}]
			\begin{axis}[height=\plotheightd\textheight,width=\plotwidthd\textwidth,xmax=0,xmin=-50,legend columns=1,legend style={draw=\legendboxdraw,fill=\legendboxfill,at={(0,0.95)},anchor=north west},legend cell align={left},
			xlabel={$z$}, ylabel={},label style={font=\normalsize},tick label style={font=\normalsize},legend image post style={scale=\legendmarkscale},legend style={nodes={scale=\legendfontscale, transform shape}},grid=none]
			\addplot[color=colorone,line width=\plotlinewidth pt,mark=\markone,mark repeat=120,mark phase=1] table [x=z,y=c,col sep=comma] 
			{data/text/c_0_s_5.csv};\addlegendentry{$c_s(z)=d_s(z),\, \varepsilon=0$};
%			\addplot[color=colortwo,line width=\plotlinewidth pt,mark=\marktwo,mark repeat=120,mark phase=41] table [x=z,y=c,col sep=comma] 
%			{data/text/c_0.01_s_5.csv};\addlegendentry{$\varepsilon=0.01$};
			\addplot[color=colortwo,line width=\plotlinewidth pt,mark=\marktwo,mark repeat=120,mark phase=81] table [x=z,y=c,col sep=comma] 
			{data/text/c_0.05_s_5.csv};\addlegendentry{$c_s(z),\, \varepsilon=0.05$};
			\addplot[color=colorthree,line width=\plotlinewidth pt,mark=\markthree,mark repeat=120,mark phase=41] table [x=z,y=d,col sep=comma] 
			{data/text/c_0.05_s_5.csv};\addlegendentry{$d_s(z),\, \varepsilon=0.05$};
			\end{axis}
		\end{tikzpicture}
		\caption{Displaying $c_s(z)$ and $d_s(z)$ as a function of $z$, for fixed $s=5$ and $\varepsilon=0$ or $\varepsilon=0.05$.}
		\label{fig:cz}
	\end{subfigure}\hfill
	\begin{subfigure}[t]{\subfigsize\textwidth}
	\centering
	\begin{tikzpicture}[every pin/.style={fill=white}]
	\begin{axis}[height=\plotheightd\textheight,width=\plotwidthd\textwidth,xmax=50,xmin=1,legend columns=1,legend style={draw=\legendboxdraw,fill=\legendboxfill,at={(1,0.35)},anchor=south east},legend cell align={left},
	xlabel={$s$}, ylabel={},label style={font=\normalsize},tick label style={font=\normalsize},legend image post style={scale=\legendmarkscale},legend style={nodes={scale=\legendfontscale, transform shape}},grid=none]
	\addplot[color=colorone,line width=\plotlinewidth pt,mark=\markone,mark repeat=10,mark phase=1] table [x=s,y=c,col sep=comma]
	{data/text/c_0.05_z_-2.csv};\addlegendentry{$z=-2$};
	\addplot[color=colortwo,line width=\plotlinewidth pt,mark=\marktwo,mark repeat=10,mark phase=1] table [x=s,y=c,col sep=comma] 
	{data/text/c_0.05_z_-20.csv};\addlegendentry{$z=-20$};
	\addplot[color=colorthree,line width=\plotlinewidth pt,mark=\markthree,mark repeat=10,mark phase=1] table [x=s,y=c,col sep=comma] 
	{data/text/c_0.05_z_-200.csv};\addlegendentry{$z=-200$};
%	\addplot[color=colorfour,line width=\plotlinewidth pt,mark=\markfour,mark repeat=9,mark phase=7] table [x=s,y=cs,col sep=comma] 
%	{data/text/c_0.05_z_-2000.csv};\addlegendentry{$z=-2000$};
	\end{axis}
	\end{tikzpicture}
	\caption{Displaying $c_s(z)$ as a function of $s$, for fixed $z=-2,-20,-200$ and $\varepsilon=0.05$.}
	\label{fig:cs}
\end{subfigure}
\end{center}
\caption{Illustration of $c_s(z)$ \cref{eq:defc} and $d_s(z)$ \cref{eq:defds} for different values of $\varepsilon$, $s$ and $z$.}
\label{fig:ceps}
\end{figure}
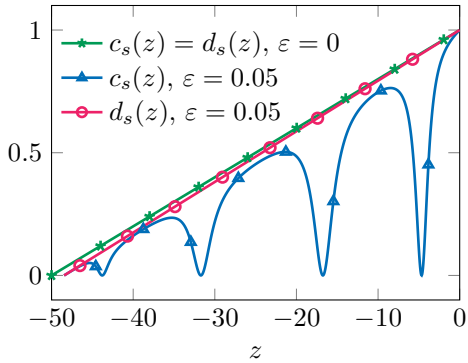
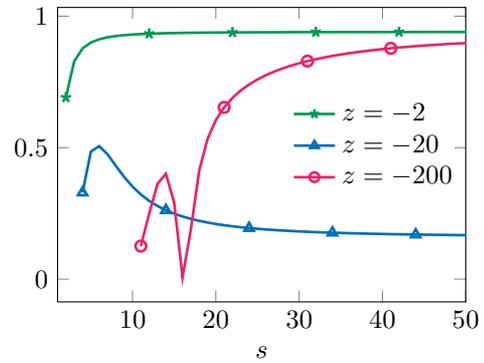

\paragraph{The PSK-$\tau$-ROCK method.}
In \cite{AAV18}, the authors design a postprocessor for the SK-ROCK scheme when applied to nonlinear Brownian dynamics. Such postprocessed scheme, called PSK-ROCK, captures the invariant measure of the ergodic process with second-order accuracy and is exact for linear equations as the Ornstein--Uhlenbeck problem. The accuracy analysis of the PSK-ROCK method is based on ergodic properties of numerical schemes for SDEs driven by diffusion processes, see for instance \cite{AVZ14,Tal90,TaT90,Vil15}. Therefore, the postprocessing strategy for the SK-ROCK scheme cannot be extended straightaway to SDEs driven by Poisson noise, as those considered here. We therefore analyze in this section the PSK-$\tau$-ROCK scheme. A full analysis for the linear problem \cref{eq:reviso} is provided, while we verify numerically in \cref{sec:numexp} that the postprocessing techniques proposed in \cref{sec:postalgo} also successfully apply to nonlinear problems.

The next theorem provides an expression for the mean and variance of the postprocessed step \cref{eq:postproc}, for a general parameter $\alpha$.
\begin{theorem}\label{thm:post}
	If $\tau\rho\leq \beta s^2$ and $\varepsilon>0$ the mean and variance of the PSK-$\tau$-ROCK scheme \cref{eq:sktaurock,eq:postproc} applied to \cref{eq:reviso} are absolutely stable and satisfy 
	\begin{equation}\label{eq:sktautockstatpost}
	\lim_{n\to\infty} \exp(\overline X_n) = \exp(X_\infty), \qquad\qquad \lim_{n\to\infty} \var(\overline X_n)  =\bar c_s(z) \var(X_\infty),
	\end{equation}
	with
	\begin{equation}\label{eq:defcbar}
	 \bar c_s(z) = c_s(z) -2\alpha^2 z
	\end{equation}
	and $c_s(z)$ as in \cref{eq:defc}.
\end{theorem}
\begin{proof}
	As $\exp(Q(X_n,\tau))=0$, from \cref{eq:postproc} follows $\lim_{n\to\infty}\exp(\overline X_n)=\lim_{n\to\infty}\exp( X_n)=\exp(X_\infty)$. From the law of total variation, for the variance we obtain
	\begin{equation}\label{eq:stabsktaurocktmp4}
		\var(\overline X_n) =\var(X_n)+\alpha^2\var(Q(X_n,\tau)). %&=\var(X_n)+\alpha_s(\varepsilon)^2(c_1-c_2)A_s(z)^n\exp(\tX_0)-2\alpha_s(\varepsilon)^2\frac{c_1c_2}{\lambda}X^T\tau
	\end{equation}
From \cref{eq:invmes,eq:stabsktaurocktmp3}, using $|A_s(z)|<1$, we deduce $\lim_{n\to\infty} \var(Q(X_n,\tau))=-2z\var(X_\infty)$. 
We conclude taking the limit in \cref{eq:stabsktaurocktmp4} and using  \cref{eq:sktautockstat}.
\end{proof}

We discuss here our choice for $\alpha$ in \cref{eq:defalpha}, which aims at providing an amplifying factor $\bar c_s(z)\leq 1$ as close to $1$ as possible.
From \cref{eq:defcbar} we note that it is possible to obtain $\bar c_s(z)=1$ for all $z\in (-\ell_s^\varepsilon,0]$ only if $c_s(z)$ is an affine function of the form $d_s(z)=1+2\alpha^2 z$. From \cref{fig:cz} we deduce that this is untrue, unless $\varepsilon=0$. However, even though $c_s(z)$ oscillates, we observe in \cref{fig:cz} that it often approaches the affine function $d_s(z)$ passing through the end points $(0,c_s(0))$ and $(-\ell_s^\varepsilon,c_s(-\ell_s^\varepsilon))$, i.e.
\begin{equation}\label{eq:defds}
	d_s(z) = c_s(0)+\frac{c_s(-\ell_s^\varepsilon)-c_s(0)}{-\ell_s^\varepsilon}z=1+\frac{1}{\ell_s^\varepsilon}z =1+\frac{\ou}{2\oz}z,
\end{equation}
where we used $c_s(0)=1$, $c_s(-\ell_s^\varepsilon)=0$\footnote{We would not have $c_s(-\ell_s^\varepsilon)=0$ if $\nu_1$ was defined exactly as for the SK-ROCK scheme \cite{AAV18}.} and $\ell_s^\varepsilon=2\oz/\ou$. Therefore, we choose $\alpha$ so that
\begin{equation}
	d_s(z)-2\alpha^2 z=1,\quad\mbox{i.e.}\quad \alpha=\frac{1}{2}\sqrt{\frac{\ou}{\oz}}.
\end{equation}
The resulting amplifying factor $\bar c_s(z)$ \cref{eq:defcbar} for the postprocessed scheme PSK-$\tau$-ROCK is displayed in \cref{fig:bcz} for different values of $\varepsilon$ and $s=5$. In \cref{fig:bcs} we display $\bar c_s(z)$ as a function of $s$ for fixed $z=-200,-20,-2$ and $\varepsilon=0.05$. We see that the variance is not amplified with this definition of $\bar c_s(z)$. Furthermore its damping remains bounded. This is in sharp contrast with the $\tau$-ROCK or Rev-$\tau$-ROCK methods (or some standard explicit or implicit $\tau$-leap methods), where the variance is either amplified or strongly damped, respectively. 
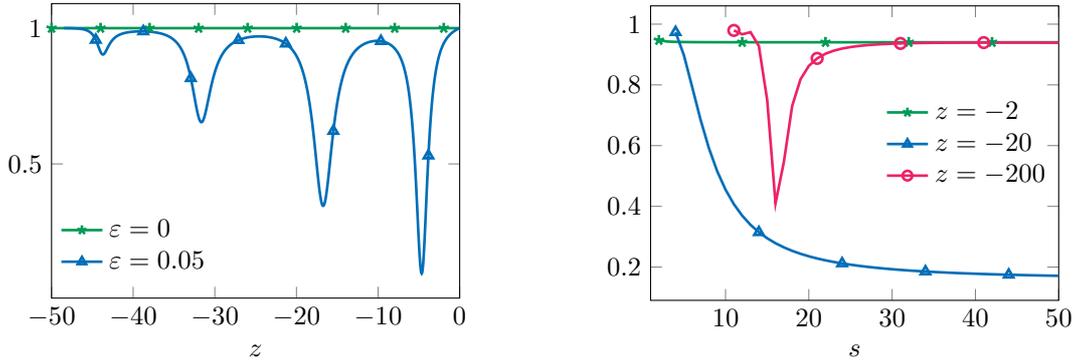
\begin{figure}
	\begin{center}
		\begin{subfigure}[t]{\subfigsize\textwidth}
			\centering
			\begin{tikzpicture}[every pin/.style={fill=white}]
			\begin{axis}[height=\plotheightd\textheight,width=\plotwidthd\textwidth,xmax=0,xmin=-50,legend columns=1,legend style={draw=\legendboxdraw,fill=\legendboxfill,at={(0,0.05)},anchor=south west},legend cell align={left},
			xlabel={$z$}, ylabel={},label style={font=\normalsize},tick label style={font=\normalsize},legend image post style={scale=\legendmarkscale},legend style={nodes={scale=\legendfontscale, transform shape}},grid=none]
			\addplot[color=colorone,line width=\plotlinewidth pt,mark=\markone, mark repeat=120,mark phase=1] table [x=z,y=bc,col sep=comma] 
			{data/text/bc_0_s_5.csv};\addlegendentry{$\varepsilon=0$};
%			\addplot[color=colortwo,line width=\plotlinewidth pt,mark=\marktwo, mark repeat=120,mark phase=41] table [x=z,y=bc,col sep=comma] 
%			{data/text/bc_0.01_s_5.csv};\addlegendentry{$\varepsilon=0.01$};
			\addplot[color=colortwo,line width=\plotlinewidth pt,mark=\marktwo, mark repeat=120,mark phase=81] table [x=z,y=bc,col sep=comma] 
			{data/text/bc_0.05_s_5.csv};\addlegendentry{$\varepsilon=0.05$};
%			\coordinate (pta) at (axis cs:-48,0.95);
%			\coordinate (ptb) at (axis cs:-43,0.75);
			\end{axis}
%			\draw[] (pta)node[circle,anchor=south,draw]{}--(ptb)node[anchor=north,draw,circle,fill=white,inner sep=0,outer sep=0]{
%				\begin{tikzpicture}[scale=\zoomsize,trim axis left,trim axis right]
%				\begin{axis}[tiny,hide axis,xlabel={},ylabel={},ticks=none,xmin=-49,xmax=-46.5]
%				\addplot[color=colorone,line width=\zoomlinewidth pt,mark=\markone,mark size=\plotmarksized,mark repeat=12,mark phase=1] table [x=z,y=bc,col sep=comma] {data/text/bc_0_s_5.csv};
%				\addplot[color=colortwo,line width=\zoomdashedlinewidth pt,mark=\marktwo,mark size=\plotmarksized,mark repeat=12,mark phase=5] table [x=z,y=bc,col sep=comma] {data/text/bc_0.01_s_5.csv};
%				\addplot[color=colorthree,line width=\zoomdashedlinewidth pt,mark=\markthree,mark size=\plotmarksized,mark repeat=12,mark phase=9] table [x=z,y=bc,col sep=comma] {data/text/bc_0.05_s_5.csv};
%				\end{axis}
%				\end{tikzpicture}
%			};
		\end{tikzpicture}
		\caption{Displaying $\bar c_s(z)$ as a function of $z$, for fixed $s=5$ and $\varepsilon=0$ or $\varepsilon=0.05$.}
		\label{fig:bcz}
	\end{subfigure} \hfill
	\begin{subfigure}[t]{\subfigsize\textwidth}
		\centering
		\begin{tikzpicture}[every pin/.style={fill=white}]
		\begin{axis}[height=\plotheightd\textheight,width=\plotwidthd\textwidth,xmax=50,xmin=1,legend columns=1,legend style={draw=\legendboxdraw,fill=\legendboxfill,at={(1,0.35)},anchor=south east},legend cell align={left},
		xlabel={$s$}, ylabel={},label style={font=\normalsize},tick label style={font=\normalsize},legend image post style={scale=\legendmarkscale},legend style={nodes={scale=\legendfontscale, transform shape}},grid=none]
		\addplot[color=colorone,line width=\plotlinewidth pt,mark=\markone,mark repeat=10,mark phase=1] table [x=s,y=bc,col sep=comma]
		{data/text/bc_0.05_z_-2.csv};\addlegendentry{$z=-2$};
		\addplot[color=colortwo,line width=\plotlinewidth pt,mark=\marktwo,mark repeat=10,mark phase=1] table [x=s,y=bc,col sep=comma] 
		{data/text/bc_0.05_z_-20.csv};\addlegendentry{$z=-20$};
		\addplot[color=colorthree,line width=\plotlinewidth pt,mark=\markthree,mark repeat=10,mark phase=1] table [x=s,y=bc,col sep=comma] 
		{data/text/bc_0.05_z_-200.csv};\addlegendentry{$z=-200$};
		%	\addplot[color=colorfour,line width=\plotlinewidth pt,mark=\markfour,mark repeat=9,mark phase=7] table [x=s,y=cs,col sep=comma] 
		%	{data/text/c_0.05_z_-2000.csv};\addlegendentry{$z=-2000$};
		\end{axis}
		\end{tikzpicture}
		\caption{Displaying $\bar c_s(z)$ as a function of $s$, for fixed $z=-2,-20,-200$ and $\varepsilon=0.05$.}
		\label{fig:bcs}
	\end{subfigure}
\end{center}
\caption{Illustration of $\bar c_s(z)$ for different values of $\varepsilon$, $s$ and $z$.}
\label{fig:bceps}
\end{figure}

We conclude this section showing that for $\varepsilon=0$ if \cref{eq:sktautockstatpost} holds then the PSK-$\tau$-ROCK scheme captures the variance exactly, that is $\bar c_s(z)=1$ for all $z$, as we see in \cref{fig:bcz}. We recall, however, that $\varepsilon=0$ may lead to instabilities.
\begin{corollary}\label{cor:alpha0}
	For $\varepsilon=0$ it holds $\bar c_s(z)=1$ for all $z\in [-\beta s^2,0]$.
\end{corollary}
\begin{proof}
	For $\varepsilon=0$ it holds $\oz=1$, $\ou=1/s^2$ and $\alpha=1/(2s)$. Moreover,
	\begin{equation}
		A_s(z)=T_s\left(1+\frac{z}{s^2}\right), \qquad\qquad B_s(z)=\frac{1}{s}U_{s-1}\left(1+\frac{z}{s^2}\right)\left(1+\frac{z}{2s^2}\right).
	\end{equation}
	The identity $(1-x^2)U_{s-1}(x)^2=1-T_s(x)^2$ implies $c_s(z)=1+z/(2 s^2)$, thus $\bar c_s(z)=1$.
\end{proof}

\IfStandalone
{
	\bibliographystyle{siam}
	\bibliography{../../../../../../LaTeX/library}
}{}

\section{Numerical experiments}\label{sec:numexp}
In this section we consider different numerical experiments in order to verify the stability and accuracy properties of the PSK-$\tau$-ROCK scheme of \cref{sec:psktrock} and also asses its efficiency compared to other $\tau$-leap methods for stiff problems. To do so, we first evaluate its accuracy on a nonstiff bistable problem, where we investigate the effect of the postprocessing step \cref{eq:postproc} on nonlinear problems. Then, we consider a mean-square stable problem containing fast variables. Next, we tackle a nonlinear problem where the fast variables have very large variance and thus the equation is non mean-square stable. Finally we consider an application to a more involved problem, namely a genetic feedback loop. In \cref{exp:MM,exp:nlrevreac,exp:genfeed} we compare the efficiency of the new scheme against classical state of the art $\tau$-leap methods for chemical kinetics.

\paragraph{Numerical methods and implementation details.}
Let us provide some general implementation details concerning the next experiments. In what follows, we often compare the PSK-$\tau$-ROCK method \cref{eq:sktaurock,eq:postproc} (also \cref{algo:pseudocode} in \cref{sec:pseudocode}) against other $\tau$-leap schemes: the SK-$\tau$-ROCK scheme defined by \cref{eq:sktaurock} but without postprocessing \cref{eq:postproc}, the explicit stabilized $\tau$-ROCK and Rev-$\tau$-ROCK method \cite{AHL10}, the implicit $\tau$-leap (Imp-$\tau$-leap) method and its postprocessed version (PImp-$\tau$-leap) \cite{RPC03}, the trapezoidal $\tau$-leap (Trap-$\tau$-leap) method \cite{CPR04} and the more recent split-step implicit $\tau$-leap (SSI-$\tau$-leap) method \cite{HMT17}. 
All the aforementioned methods are implemented in C++ using the \texttt{Eigen} library \cite{GuB10} for the linear algebra routines. 
The recently introduced split step method \cite{RKV19} is not considered in the following experiments; nevertheless, notice that its cost is roughly twice the cost of the Imp-$\tau$-leap method.

For the PSK-$\tau$-ROCK and SK-$\tau$-ROCK method we always use $\varepsilon=0.05$ as damping parameter, i.e. the standard choice for explicit stabilized methods, and $\alpha$ as in \cref{eq:defalpha}. The number of stages $s$ is chosen before each step according to the condition $\tau\rho\leq \beta s^2$ with $\beta=2-4\varepsilon/3$, where the spectral radius $\rho$ of the Jacobian of $f$ (see \cref{eq:deffQtauleap}) is approximated employing a nonlinear power method \cite{Lin72,Ver80}. If relevant, we report the number of iterations of this nonlinear power method (PM). Here the condition $\tau\rho\leq \beta s^2$ guarantees stability and as $\varepsilon$ is small $\beta s^2$ is a good approximation of the true stability domains' size $\ell_s^\varepsilon=2\oz/\ou$.

For the $\tau$-ROCK and Rev-$\tau$-ROCK method the damping parameter $\varepsilon=\varepsilon(s)$ depends on $s$ and might be large, therefore $\beta s^2$ does not approximate the stability domain's size accurately. Hence, for these methods, the number of stages $s$ is the smallest integer verifying $\tau\rho\leq \ell_s^\varepsilon$, where $\ell_s^\varepsilon=2\oz/\ou$ is the exact size of the stability domain of a scheme with $s$ stages and damping $\varepsilon$. The damping parameter is chosen according to the strategy described in \cite{AbL08} for the $\tau$-ROCK method and according to \cite[eq. 4.32]{AHL10} for the Rev-$\tau$-ROCK method. 

Finally, for the postprocessing step of the PImp-$\tau$-leap scheme we use ten steps of size $\delta\tau=0.2/\rho$ of the Imp-$\tau$-leap method (note that the relaxation time of the fast variables is proportional to $1/\rho$).

We conclude this paragraph commenting on negative populations. As the Poisson random variables are unbounded, it is a common issue of $\tau$-leap methods that negative populations arise when too many reactions fire during one step, compared to the available number of reactants. In the literature, numerous strategies have been developed in order to guarantee positive populations \cite{CGP05c,CVK05,RaS07,TiB04,YRS11}; however, developing or adapting such a strategy for the SK-$\tau$-ROCK methods is not the focus of this paper. Therefore, in the following numerical experiments we employ the usual trick of considering the absolute value of the components of $\bx$ whenever $Q(\bx,\tau)$ is evaluated.

\subsection{Accuracy experiment on the nonstiff Schlögl model}
In this first numerical experiment we want to investigate the effect of the postprocessing procedure \cref{eq:postproc} on the accuracy of the solution. To do so, we compare the PSK-$\tau$-ROCK method with the non postprocessed SK-$\tau$-ROCK method on the nonstiff nonlinear Schlögl model:
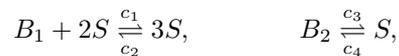
\begin{equation}\label{eq:schlogl}
B_1+2S \overset{c_1}{\underset{c_2}{\rightleftharpoons}} 3 S,\qquad\qquad B_2 \overset{c_3}{\underset{c_4}{\rightleftharpoons}} S,
\end{equation}
where $B_1$ and $B_2$ are buffered species whose populations are kept constant at $N_1=10^5$ and $N_2=2\cdot 10^5$ molecules, respectively. The state vector $X(t)$ represents the number of molecules of $S$, we set $X(0)=250$ and the final time $T=50$. The state-change vectors are $\bnu_1=\bnu_3=1$, $\bnu_2=\bnu_4=-1$, the propensity functions are
\begin{equation}
	a_1(x)=c_1 N_1x(x-1)/2,\qquad a_2(x)=c_2 x(x-1)(x-2)/6,\qquad a_3(x)=c_3 N_2,\qquad a_4(x)=c_4 x
\end{equation}
and we set $c_1=3\cdot 10^{-7}$, $c_2=10^{-4}$, $c_3=10^{-3}$ and $c_4=3.5$. For this choice of parameters the solution $X$ has a bistable distribution. 

In this experiment we fix $\tau=0.5$ and as the problem is nonstiff the PSK-$\tau$-ROCK and SK-$\tau$-ROCK methods are stable already for $s=1$; however, note that even though $s=1$ the PSK-$\tau$-ROCK and SK-$\tau$-ROCK methods do not correspond to the standard $\tau$-leap scheme \cref{eq:tauleap}. Using $10^6$ samples, we estimate the probability density function (pdf) of $X(T)$ approximated by the PSK-$\tau$-ROCK or the SK-$\tau$-ROCK method with $s=1$. We display the results in \cref{fig:schlogl_pdf} and compare them against a reference pdf computed with the SSA. We observe that the PSK-$\tau$-ROCK method matches very well the reference solution, while the SK-$\tau$-ROCK scheme tends to cluster its solutions too close to the two stable points; both results are in line with the results of \cref{sec:analysis}.

Now we want to investigate the accuracy of the schemes with respect to the number of stages $s$. To do so, we define the density distance area (DDA) \cite{CaP06}
\begin{equation}
\operatorname{DDA} = \int_{\Rb}|\overline p(x)-p(x)|\dif x,
\end{equation}
where $p(x)$ is the probability density function of $X(T)$ computed by the SSA and $\overline p(x)$ the one computed by the PSK-$\tau$-ROCK or SK-$\tau$-ROCK method. We display in \cref{fig:schlogl_dda} the DDA of the two schemes for different stages $s$ and fixed $\tau=0.5$, where the pdf $\overline p(x)$ is estimated using $10^6$ samples. We note that the SK-$\tau$-ROCK method becomes more accurate as $s$ increases. In contrast, the PSK-$\tau$-ROCK method is accurate already for low values of $s$.
%; this result is in line with \cref{fig:cs}, where we see that the amplification factor $c_s(z)$ for $z=-2$ gets closer to the optimal value $1$ as $s$ increases. In contrast, the PSK-$\tau$-ROCK method's accuracy does not vary much with $s$; this is also in line with \cref{fig:bcs}, where we note that the amplification factor $\bar c_s(z)$ for $z=-2$ remains close to one even for larger $s$. 
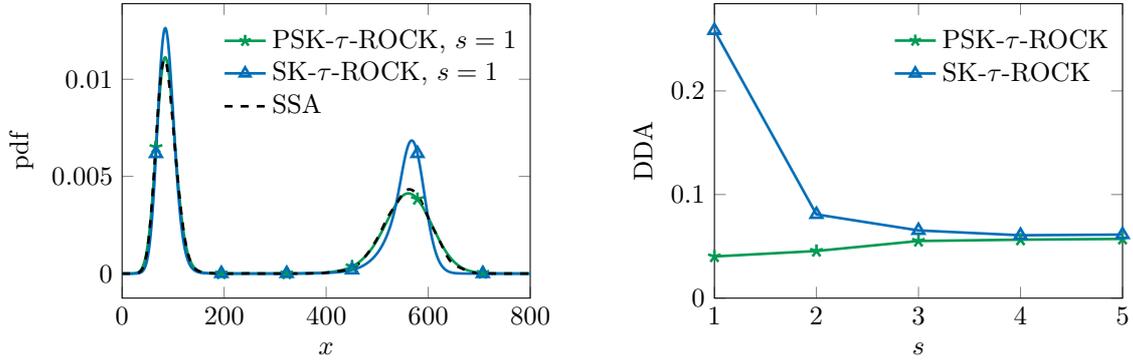
\begin{figure}
	\begin{center}
		\begin{subfigure}[t]{\subfigsize\textwidth}
			\centering
			\begin{tikzpicture}
			\begin{axis}[height=\plotheightd\textheight,width=\plotwidthd\textwidth,legend columns=1,legend style={draw=\legendboxdraw,fill=\legendboxfill,at={(1,0.95)},anchor=north east},log basis x={2},log basis y={2},legend cell align={left},xmin=0,xmax=800,
			xlabel={$x$}, ylabel={pdf},label style={font=\normalsize},tick label style={font=\normalsize},y tick label style={/pgf/number format/fixed,/pgf/number format/precision=4},legend image post style={scale=\legendmarkscale},legend style={nodes={scale=\legendfontscale, transform shape}},scaled ticks=false]
			\addplot[color=colorone,solid,line width=\plotlinewidth pt,mark=\markone,mark size=\plotmarksizeu pt, mark repeat=160,mark phase=84] table [x=res1,y=res2,col sep=comma] 
			{data/Schlogl/distribution_pstr_mc_1e6_dt_0p5.csv};\addlegendentry{PSK-$\tau$-ROCK, $s=1$}
			\addplot[color=colortwo,solid,line width=\plotlinewidth pt,mark=\marktwo,mark size=\plotmarksizeu pt, mark repeat=160,mark phase=84] table [x=res1,y=res2,col sep=comma] 
			{data/Schlogl/distribution_str_mc_1e6_dt_0p5.csv};\addlegendentry{SK-$\tau$-ROCK, $s=1$}
			\addplot[color=black,dashed,line width=\plotlinewidth pt,mark=none,mark size=\plotmarksizeu pt, mark repeat=160,mark phase=84] table [x=res1,y=res2,col sep=comma] 
			{data/Schlogl/distribution_ssa_mc_1e8.csv};\addlegendentry{SSA}
			\end{axis}
			\end{tikzpicture}
			\caption{Comparing the reference pdf computed by the SSA against the pdf computed by the two methods, with $s=1$.}
			\label{fig:schlogl_pdf}
		\end{subfigure}\hfill
	\begin{subfigure}[t]{\subfigsize\textwidth}
		\centering
		\begin{tikzpicture}
		\begin{axis}[height=\plotheightd\textheight,width=\plotwidthd\textwidth,legend columns=1,legend style={draw=\legendboxdraw,fill=\legendboxfill,at={(1,0.95)},anchor=north east},log basis x={2},log basis y={2},legend cell align={left},xmin=1,xmax=5,ymin=0,
		xlabel={$s$}, ylabel={DDA},label style={font=\normalsize},tick label style={font=\normalsize},legend image post style={scale=\legendmarkscale},legend style={nodes={scale=\legendfontscale, transform shape}},]
		\addplot[color=colorone,solid,line width=\plotlinewidth pt,mark=\markone,mark size=\plotmarksizeu pt] table [x=s,y=pstr,col sep=comma] 
		{data/Schlogl/s_VS_dda.csv};\addlegendentry{PSK-$\tau$-ROCK}
		\addplot[color=colortwo,solid,line width=\plotlinewidth pt,mark=\marktwo,mark size=\plotmarksizeu pt] table [x=s,y=str,col sep=comma] 
		{data/Schlogl/s_VS_dda.csv};\addlegendentry{SK-$\tau$-ROCK}
		\end{axis}
		\end{tikzpicture}
		\caption{Errors committed on the pdf by the two methods, as functions of the stages $s$.}
		\label{fig:schlogl_dda}
	\end{subfigure}
	\end{center}
	\caption{Schlögl reaction. Approximation of the probability density function (pdf) of $X(T)$ computed by the PSK-$\tau$-ROCK and SK-$\tau$-ROCK methods, with $\tau=0.5$ and different stages $s$.}
	\label{fig:conv}
\end{figure}

\subsection{Efficiency experiment on the Michaelis--Menten system}\label{exp:MM}
Here we consider the Michaelis--Menten system already described in \cref{ex:MM}, where we set $\bX(0)=(3000,120,0,0)^\top$ and the reaction rate constants as $c_1=1.66\cdot 10^{-3}$, $c_2=10^{-4}$ and $c_3=10^3$. With this set of coefficients, the variables $X_1,X_4$ are slow and $X_2,X_3$ are fast; however, this is a mean-square stable problem and therefore all variances tend to zero. A typical solution of the Michaelis--Menten system is displayed in \cref{fig:MMsol}. For this model, the quantities of interest are the slow variables $X_1,X_4$; therefore the byproducts $X_2,X_3$ are neglected in what follows.
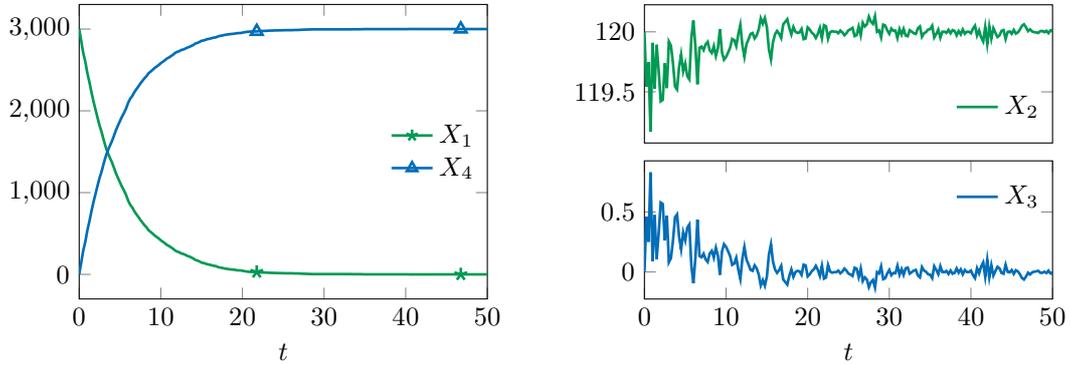
\begin{figure}
	\begin{center}
		\begin{subfigure}[b]{\subfigsize\textwidth}
			\centering
			\begin{tikzpicture}
				\begin{axis}[height=\plotheightd\textheight,width=\plotwidthd\textwidth,legend columns=1,legend style={draw=\legendboxdraw,fill=\legendboxfill,at={(1,0.5)},anchor=east},legend cell align={left},xmin=0,xmax=50,
					xlabel={$t$}, ylabel={},label style={font=\normalsize},tick label style={font=\normalsize},y tick label style={/pgf/number format/fixed,/pgf/number format/precision=4},legend image post style={scale=\legendmarkscale},legend style={nodes={scale=\legendfontscale, transform shape}},scaled ticks=false]
					\addplot[color=colorone,solid,line width=\plotlinewidth pt,mark=\markone,mark size=\plotmarksizeu pt, mark repeat=100,mark phase=88] table [x=t,y=x1,col sep=comma] 
					{data/MichaelisMenten/MichaelisMenten_path.csv};\addlegendentry{$X_1$}
					\addplot[color=colortwo,solid,line width=\plotlinewidth pt,mark=\marktwo,mark size=\plotmarksizeu pt, mark repeat=100,mark phase=88] table [x=t,y=x4,col sep=comma] 
					{data/MichaelisMenten/MichaelisMenten_path.csv};\addlegendentry{$X_4$}
				\end{axis}
			\end{tikzpicture}
%		\caption{Michaelis--Menten. Illustration of one solution's sample.}
%		\label{fig:MMsol_14}
		\end{subfigure}
	\begin{subfigure}[b]{\subfigsize\textwidth}
		\centering
		\begin{tikzpicture}
			\begin{axis}[height=\plotheightdhalf\textheight,width=\plotwidthd\textwidth,at={(0,130)},legend columns=1,legend style={draw=\legendboxdraw,fill=\legendboxfill,at={(1,0.1)},anchor=south east},legend cell align={left},xmin=0,xmax=50,xmajorticks=false,%ymin=119,ymax=120.5,ytick={119,120},
				xlabel={}, ylabel={},label style={font=\normalsize},tick label style={font=\normalsize},y tick label style={/pgf/number format/fixed,/pgf/number format/precision=4},legend image post style={scale=\legendmarkscale},legend style={nodes={scale=\legendfontscale, transform shape}},scaled ticks=false]
				\addplot[color=colorone,solid,line width=\plotlinewidth pt,mark=none,mark size=\plotmarksizeu pt, mark repeat=100,mark phase=88] table [x=t,y=x2,col sep=comma] 
				{data/MichaelisMenten/MichaelisMenten_path.csv};\addlegendentry{$X_2$}
			\end{axis};
%		\end{tikzpicture} \\
%	\begin{tikzpicture}
		\begin{axis}[height=\plotheightdhalf\textheight,width=\plotwidthd\textwidth,at={(0,0)},legend columns=1,legend style={draw=\legendboxdraw,fill=\legendboxfill,at={(1,0.9)},anchor=north east},legend cell align={left},xmin=0,xmax=50,
			xlabel={$t$}, ylabel={},label style={font=\normalsize},tick label style={font=\normalsize},y tick label style={/pgf/number format/fixed,/pgf/number format/precision=4},legend image post style={scale=\legendmarkscale},legend style={nodes={scale=\legendfontscale, transform shape}},scaled ticks=false]
			\addplot[color=colortwo,solid,line width=\plotlinewidth pt,mark=none,mark size=\plotmarksizeu pt, mark repeat=100,mark phase=88] table [x=t,y=x3,col sep=comma] 
			{data/MichaelisMenten/MichaelisMenten_path.csv};\addlegendentry{$X_3$}
		\end{axis}
	\end{tikzpicture}
%		\caption{Michaelis--Menten. Illustration of one solution's sample.}
%		\label{fig:MMsol_23}
	\end{subfigure}
	\end{center}
	\caption{Michaelis--Menten. Illustration of a typical solution.}
	\label{fig:MMsol}
\end{figure}

In this experiment, we want to compare the accuracy and efficiency of the PSK-$\tau$-ROCK method against the other $\tau$-leap methods listed at the beginning of \cref{sec:numexp}. To do so, we fix $\tau=0.25$, integrate the equations with the different $\tau$-leap methods and compare the expectations and standard deviations, computed over $10^6$ samples, against a reference solution computed with the SSA. Also, we verify accuracy at the transient time $T=5$ and at the equilibrium time $T=50$.
The results are reported in \cref{tab:MMstat}.
We observe that all the schemes approximate relatively well the slow variables $X_1,X_4$, as the step size $\tau=0.25$ is small enough to resolve them and the fast variables have too small variance to perturb the accuracy of slow dynamics. Comparing the PSK-$\tau$-ROCK and SK-$\tau$-ROCK schemes, we observe that for this mean-square stable problem the postprocessing procedure \cref{eq:postproc} has no effect on the solution's accuracy.
\begin{table}
	\centering
	\begin{subtable}{\textwidth}
		\centering
		\begin{tabular}{lcccc}
			\toprule
			& $X_1(5)$ & $X_4(5)$ & $X_1(50)$ & $X_4(50)$ \\
			\midrule
			SSA (reference) &1111.6 & 1888.2 & 0.14 & 2999.86 \\
			PSK-$\tau$-ROCK  & 1093.2 & 1906.6  & 0.12 & 2999.88 \\
			SK-$\tau$-ROCK  & 1093.2 & 1906.6 & 0.12 & 2999.88 \\ 
			$\tau$-ROCK & 1104.8 & 1895.0 & 0.13 & 2999.85  \\
			Rev-$\tau$-ROCK  & 1104.3 & 1895.5  & 0.13 & 2999.87  \\
			Imp-$\tau$-leap & 1138.3 & 1861.4 & 0.02 & 2999.98 \\
			PImp-$\tau$-leap & 1138.3  & 1861.6  & 0.02 & 2999.98  \\
			Trap-$\tau$-leap & 1111.2 & 1888.7  & -0.01 & 3000.02  \\
			SSI-$\tau$-leap  & 1138.4  & 1861.4 & 0.17 & 3000.29  
		\end{tabular}
		\caption{Empirical means.}
		\label{tab:MMmean}
	\end{subtable}\\ \vspace{0.2cm}
	\begin{subtable}{\textwidth}
		\centering
		\begin{tabular}{lcccc}
			\toprule
			& $X_1(5)$ & $X_4(5)$  & $X_1(50)$ & $X_4(50)$  \\
			\midrule
			SSA (reference) & 26.4 & 26.4  & 0.38 & 0.38  \\
			PSK-$\tau$-ROCK  & 26.6  & 26.7  & 0.44 & 0.45  \\
			SK-$\tau$-ROCK & 26.6 & 26.7  & 0.44 & 0.44 \\
			$\tau$-ROCK & 26.7 & 27.4 & 0.45 & 6.14  \\
			Rev-$\tau$-ROCK & 25.9  & 25.9 & 0.44 & 0.44  \\
			Imp-$\tau$-leap  & 26.5  & 26.5  & 0.62 & 0.62 \\
			PImp-$\tau$-leap & 26.5  & 26.5  & 0.62 & 0.62  \\
			Trap-$\tau$-leap & 26.7 & 26.7  & 0.52 & 0.52 \\
			SSI-$\tau$-leap  & 27.2 & 27.2 & 0.47 & 0.77  
		\end{tabular}
		\caption{Empirical standard deviations.}
		\label{tab:MMstd}
	\end{subtable}
	\caption{Michaelis--Menten. Empirical means and standard deviations of $X_1,X_4$ at times $T=5$ and $T=50$.}
	\label{tab:MMstat}
\end{table}

\Cref{tab:MMtime} displays the computational times of the different methods, together with the average number of stages $s$ and the average damping parameter $\varepsilon$. Moreover, for the explicit stabilized schemes we also display the average number of iterations, per time step, needed to approximate the spectral radius $\rho$ with the nonlinear power method (PM) and for the implicit methods we display the average number of iterations needed by the Newton (N) method. The PSK-$\tau$-ROCK and SK-$\tau$-ROCK methods are the most efficient schemes. 
\begin{table}
	\centering
	\begin{tabular}{lcccr}
		\toprule
		& $s$ & $\varepsilon$ & PM/N iter & CPU [sec.] \\
		\midrule
		PSK-$\tau$-ROCK & 12 & 0.05 & 1.01 & 99 \\ 
		SK-$\tau$-ROCK & 12 & 0.05 & 1.01 & 98 \\ 
		$\tau$-ROCK & 24 & 8.55 & 1.01 & 187 \\ 
		Rev-$\tau$-ROCK & 23 & 7.33 & 1.01 & 179 \\ 
		Imp-$\tau$-leap & 1 & - & 1.59 & 250 \\ 
		PImp-$\tau$-leap & 1 & - & 1.59 & 258\\ 
		Trap-$\tau$-leap & 1 & - & 1.58 &  250 \\ 
		SSI-$\tau$-leap  & 2 & - & 1.96 & 294 
	\end{tabular}
	\caption{Michaelis--Menten. Per step average of number of stages, of damping parameter, of nonlinear power method (PM) or Newton (N) iterations and total computational times taken by the different methods.}
	\label{tab:MMtime}
\end{table}

\subsection{Stiff nonlinear reversible reaction}\label{exp:nlrevreac}
In this experiment we consider the stiff nonlinear reversible reaction
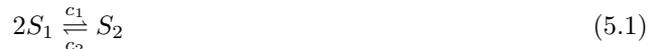
\begin{equation}\label{eq:nlrevreac}
2S_1\overset{c_1}{\underset{c_2}{\rightleftharpoons}}S_2
\end{equation}
with $c_1=50$, $c_2=10^3$ and $\bX(0)=(400,3990)^\top$. In this setting, \cref{eq:nlrevreac} is at equilibrium and we can illustrate the accuracy of the schemes in capturing the invariant measure of the system, moreover the variances are large and thus the problem is not mean-square stable. As the quantity $X_C=X_1+2X_2$ is constant over time we can eliminate $X_2$ from the system and let $X=X_1$. The propensity functions and state-change vectors are
\begin{equation}
	a_1(x)=c_1 x(x-1)/2,\qquad a_2(x)=c_2(X_C-x)/2,\qquad \bnu_1=-2, \qquad \bnu_2=2.
\end{equation}
We fix $T=0.2$, $\tau=0.01$ and integrate the system with the $\tau$-leap methods listed at the beginning of the section; the results are reported in \cref{tab:nlrevreac}, where we use $10^6$ samples to approximate the statistics. 
\begin{table}
	\centering
	\begin{tabular}{lrrrcr}
		\toprule
		& $\exp(X)$ & Std$(X)$ & $s$  &PM/N iter & CPU [sec.] \\ % & Speed-up\\
		\midrule
		SSA (reference) & 399.4 & 19.8 & - & - &  17325 \\% & 1x \\ 
		PSK-$\tau$-ROCK & 397.4 & 19.7  & 16.2 &  1.71 &  12 \\% & 1333x \\ 
		SK-$\tau$-ROCK & 397.3 & 8.6  & 16.2 &  1.71 &  12 \\% & 1333x \\
%		$\tau$-ROCK & 800 & 3500 & 396.5 & 561  & (10505)  \\% & <2x\\
%		Rev-$\tau$-ROCK & 800 & 3500 & 400 & 0  & (8973) \\% & <2x\\
		Imp-$\tau$-leap &  400.0 & 1.4  & 1 &  1.04 & 10 \\% & 1733x\\
		PImp-$\tau$-leap & 399.7 & 18.6  & 1 &  1.07 &  16 \\% & 963x\\ 
		Trap-$\tau$-leap & 399.8 & 11.23  & 1 &  1.75 &  15 \\% & 1238x\\
		SSI-$\tau$-leap & 1.3 & 708  & 2 & 4.17 &  29 \\% & 597x\\
	\end{tabular}
	\caption{Nonlinear reversible reaction. For different methods, we report the mean and standard deviation of $X$, the per step average of number of stages, the average number of nonlinear power method (PM) or Newton (N) iterations and the total computational time.}
	\label{tab:nlrevreac}
\end{table}
For the $\tau$-ROCK scheme we could not use the standard strategy to define the number of stages, as for this nonlinear problem with large variance we found that it does not guarantee stability. 
Searching for a set of parameters $\varepsilon,s$ providing at least 1\% of stable paths we found $\varepsilon=3500$, $s=800$, which is an unusable number of stages due to round-off errors \cite{VHS90}. The same holds for Rev-$\tau$-ROCK. We therefore did not include these methods in our numerical comparisons. All the other schemes were successful in $100\%$ of the Monte-Carlo iterations.
We observe in \cref{tab:nlrevreac} that the SSI-$\tau$-leap scheme is completely off. The other methods approximate well the expectation, while only PSK-$\tau$-ROCK and PImp-$\tau$-leap provide good approximations to the standard deviation, with PSK-$\tau$-ROCK being the most accurate at smaller cost. For this very small problem, the Imp-$\tau$-leap method is slightly faster than the PSK-$\tau$-ROCK scheme; in contrast, it is significantly less precise.
In \cite{CPR04} it is shown that the trapezoidal rule captures the exact invariant measure of \cref{eq:reviso}, i.e. a linear problem. This example shows that it fails in capturing the statistics of nonlinear problems as \cref{eq:nlrevreac}. The same phenomenon is observed in \cite{LAE08} for the trapezoidal rule for stochastic differential equations driven by diffusion processes.

\subsection{Numerical experiment on a genetic positive feedback loop}\label{exp:genfeed}
We consider here a stiff biological system modeling a genetic positive feedback loop. The system describes the production of a protein, which auto regulates its production rate by binding to its gene promoter \cite{RaS07,YRS11}. This system is described by the set of reactions
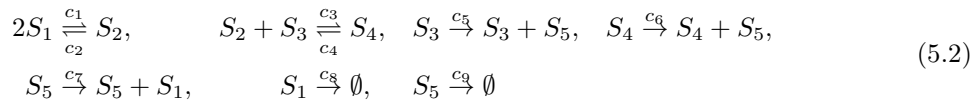
\begin{equation}\label{eq:genloop}
\begin{aligned}
2S_1&\overset{c_1}{\underset{c_2}{\rightleftharpoons}}S_2,&
S_2+S_3&\overset{c_3}{\underset{c_4}{\rightleftharpoons}}S_4,&
S_3&\overset{c_5}{\rightarrow}S_3+S_5,&
S_4&\overset{c_6}{\rightarrow}S_4+S_5,\\
S_5&\overset{c_7}{\rightarrow}S_5+S_1,&
S_1&\overset{c_8}{\rightarrow}\emptyset,&
S_5&\overset{c_9}{\rightarrow}\emptyset &&
\end{aligned}
\end{equation}
and the rate constants are given by $c_1=50$, $c_2=10^3$, $c_3=50$, $c_4=10^3$, $c_5=1$, $c_6=10$, $c_7=3$, $c_8=1$ and $c_9=6$. We set the final time to $T=100$ and the initial condition as $X(0)=(10,0,20,0,0)^\top$.

We consider the same methods as in the previous experiments and apply them to \cref{eq:genloop} using a fixed step size $\tau=0.05$ and computing $10^5$ samples. 
We start by observing that the system contains the nonlinear reversible reaction \cref{eq:nlrevreac}, which induces significant fluctuations in $X_1,X_2$. The same reaction has been considered in \cref{exp:nlrevreac}, where we found that at equilibrium the variance of $X_1,X_2$ is relatively large and the $\tau$-ROCK and Rev-$\tau$-ROCK methods need an excessively large number of stages and damping parameter in order to be stable. For \cref{eq:genloop} the same considerations hold and indeed we found that simulation of \cref{eq:genloop} with the $\tau$-ROCK and Rev-$\tau$-ROCK schemes is too often unstable; therefore, we will discard these methods from the rest of the experiment. For this example, the Trap-$\tau$-leap and the SSI-$\tau$-leap schemes had severe convergence issues in the Newton method, therefore their results are not reported neither. 

In \cref{tab:genloopmen} we note that all methods are accurate in sampling the means. In \cref{tab:genloopstdv} we observe that the PSK-$\tau$-ROCK is the most precise in estimating the variance of all the species. 
\begin{table}
	\centering
\begin{subtable}{\textwidth}
	\centering
	\begin{tabular}{l l l l l l}
		\toprule
		& $X_1$ & $X_2$ & $X_3$ & $X_4$ & $X_5$ \\
		\midrule
		SSA (reference) & 92.2 & 213 & 1.72 & 18.3 & 30.8 \\
		PSK-$\tau$-ROCK & 92.3 & 211 & 1.75 & 18.3 & 30.7 \\
		SK-$\tau$-ROCK & 92.3 & 211 & 1.75 & 18.3 & 30.8 \\ 
		Imp-$\tau$-leap & 92.1 & 210 & 1.75 & 18.3 & 30.7 \\
		PImp-$\tau$-leap & 91.9 & 210 & 1.74 & 18.3 & 30.7
	\end{tabular}
	\caption{Empirical means.}
	\label{tab:genloopmen}
\end{subtable}\\ \vspace{0.2cm}
\begin{subtable}{\textwidth}
	\centering
	\begin{tabular}{l l l l l l}
		\toprule
		& $X_1$ & $X_2$ & $X_3$ & $X_4$ & $X_5$ \\
		\midrule
		SSA (reference) & 9.87 & 18.0 & 1.26 & 1.26 & 5.55 \\
		PSK-$\tau$-ROCK  & 9.78 & 18.7 & 1.21 & 1.21 & 5.55 \\
		SK-$\tau$-ROCK & 7.24 & 18.3 & 0.72 & 0.72 & 5.54 \\
		Imp-$\tau$-leap & 4.03 & 18.3 & 0.16 & 0.16 & 5.18 \\
		PImp-$\tau$-leap & 9.41 & 18.8 & 1.21 & 1.21 & 5.19 
	\end{tabular}
	\caption{Empirical standard deviations.}
	\label{tab:genloopstdv}
\end{subtable}
	\caption{Genetic positive feedback loop. Empirical means and standard deviations.}
	\label{tab:genloopmstd}
\end{table}
The PImp-$\tau$-leap scheme seems to be almost as precise as PSK-$\tau$-ROCK, however we see in \cref{tab:genloopcpu} that it is slower. Indeed, in \cref{tab:genloopcpu} we see that PSK-$\tau$-ROCK is not only the most accurate scheme but also the fastest one, together with SK-$\tau$-ROCK. The speed-up is significantly larger than for the implicit methods.
\begin{table}
	\centering
	\begin{tabular}{lccrr}
		\toprule
		& s & PM/N iter & CPU [sec.] & Speed-up\\
		\midrule
		SSA (reference) & - & - & 49150 & 1x \\
		PSK-$\tau$-ROCK  & 21 & 1.75 &  210 &  234x \\
		SK-$\tau$-ROCK & 21 & 1.75 & 208  & 236x\\
		Imp-$\tau$-leap & 1 & 1.99 &  427 & 115x \\ 
		PImp-$\tau$-leap & 1 & 1.99 &  433 & 114x \\ 
	\end{tabular}
	\caption{Genetic positive feedback loop. For different methods, we report the per step average of the number of stages, the average number of nonlinear power method (PM) or Newton (N) iterations, the computational time and the speed-up compared to the SSA.}
	\label{tab:genloopcpu}
\end{table}

\IfStandalone
{
	\bibliographystyle{siam}
	\bibliography{../../../../../../LaTeX/library}

\begin{thebibliography}{10}
\expandafter\ifx\csname url\endcsname\relax
  \def\url#1{\texttt{#1}}\fi
\expandafter\ifx\csname doi\endcsname\relax
  \def\doi#1{\burlalt{doi:#1}{http://dx.doi.org/#1}}\fi
\expandafter\ifx\csname urlprefix\endcsname\relax\def\urlprefix{URL }\fi
\expandafter\ifx\csname href\endcsname\relax
  \def\href#1#2{#2}\fi
\expandafter\ifx\csname burlalt\endcsname\relax
  \def\burlalt#1#2{\href{#2}{#1}}\fi

\bibitem{Abd02}
A.~Abdulle.
\newblock {Fourth order Chebyshev methods with recurrence relation}.
\newblock {\em SIAM Journal on Scientific Computing}, 23(6):2041--2054, 2002.

\bibitem{AAV18}
A.~Abdulle, I.~Almuslimani, and G.~Vilmart.
\newblock {Optimal explicit stabilized integrator of weak order one for stiff
  and ergodic stochastic differential equations}.
\newblock {\em Siam Journal on Uncertainty Quantification}, 6(2):937--964,
  2018.

\bibitem{AbC07}
A.~Abdulle and S.~Cirilli.
\newblock {Stabilized methods for stiff stochastic systems}.
\newblock {\em Comptes Rendus Math{\'{e}}matique. Acad{\'{e}}mie des Sciences.
  Paris}, 345(10):593--598, 2007.

\bibitem{AbC08}
A.~Abdulle and S.~Cirilli.
\newblock {S-ROCK: Chebyshev methods for stiff stochastic differential
  equations}.
\newblock {\em SIAM Journal on Scientific Computing}, 30(2):997--1014, 2008.

\bibitem{AHL10}
A.~Abdulle, Y.~Hu, and T.~Li.
\newblock {Chebyshev Methods with Discrete Noise: the $\tau$-ROCK Methods}.
\newblock {\em Journal of Computational Mathematics}, 28(2):195--217, 2010.

\bibitem{AbL08}
A.~Abdulle and T.~Li.
\newblock {S-ROCK methods for stiff It{\^{o}} SDEs}.
\newblock {\em Communications in Mathematical Sciences}, 6(4):845--868, 2008.

\bibitem{AbM01}
A.~Abdulle and A.~A. Medovikov.
\newblock {Second order Chebyshev methods based on orthogonal polynomials}.
\newblock {\em Numerische Mathematik}, 18:1--18, 2001.

\bibitem{AVZ13b}
A.~Abdulle, G.~Vilmart, and K.~C. Zygalakis.
\newblock {Weak second order explicit stabilized methods for stiff stochastic
  differential equations}.
\newblock {\em SIAM Journal on Scientific Computing}, 35(4):A1792--A1814, 2013.

\bibitem{AVZ14}
A.~Abdulle, G.~Vilmart, and K.~C. Zygalakis.
\newblock {High order numerical approximation of the invariant measure of
  ergodic SDEs}.
\newblock {\em SIAM Journal on Numerical Analysis}, 52(4):1600--1622, 2014.

\bibitem{Blu15}
A.~Blumenthal.
\newblock {\em {Stabilized Numerical Methods for Stochastic Differential
  Equations driven by Diffusion and Jump-Diffusion Processes}}.
\newblock PhD thesis, EPFL, 2015.
\newblock \doi{10.5075/epfl-thesis-6771}.

\bibitem{But69}
J.~C. Butcher.
\newblock {The effective order of Runge-Kutta methods}.
\newblock In {\em Numerical Solution of Differential Equations}, pages
  133--139, Berlin, 1969. Springer.

\bibitem{CGP05c}
Y.~Cao, D.~T. Gillespie, and L.~R. Petzold.
\newblock {Avoiding negative populations in explicit Poisson tau-leaping}.
\newblock {\em Journal of Chemical Physics}, 123(5), 2005.

\bibitem{CGP05a}
Y.~Cao, D.~T. Gillespie, and L.~R. Petzold.
\newblock {The slow-scale stochastic simulation algorithm}.
\newblock {\em Journal of Chemical Physics}, 122(1), 2005.

\bibitem{CGP07}
Y.~Cao, D.~T. Gillespie, and L.~R. Petzold.
\newblock {Adaptive explicit-implicit tau-leaping method with automatic tau
  selection}.
\newblock {\em Journal of Chemical Physics}, 126(22), 2007.

\bibitem{CaP06}
Y.~Cao and L.~Petzold.
\newblock {Accuracy limitations and the measurement of errors in the stochastic
  simulation of chemically reacting systems}.
\newblock {\em Journal of Computational Physics}, 212(1):6--24, 2006.

\bibitem{CPR04}
Y.~Cao, L.~R. Petzold, M.~Rathinam, and D.~T. Gillespie.
\newblock {The numerical stability of leaping methods for stochastic simulation
  of chemically reacting systems}.
\newblock {\em Journal of Chemical Physics}, 121(24):12169--12178, 2004.

\bibitem{CVK05}
A.~Chatterjee, D.~G. Vlachos, and M.~A. Katsoulakis.
\newblock {Binomial distribution based $\tau$-leap accelerated stochastic
  simulation}.
\newblock {\em Journal of Chemical Physics}, 122(2), 2005.

\bibitem{ELV05b}
W.~E, D.~Liu, and E.~Vanden-Eijnden.
\newblock {Nested stochastic simulation algorithm for chemical kinetic systems
  with disparate rates}.
\newblock {\em Journal of Chemical Physics}, 123(19), 2005.

\bibitem{Gil76}
D.~T. Gillespie.
\newblock {A general method for numerically simulating the stochastic time
  evolution of coupled chemical reactions}.
\newblock {\em Journal of Computational Physics}, 22:403--434, 1976.

\bibitem{Gil77}
D.~T. Gillespie.
\newblock {Exact stochastic simulation of coupled chemical reactions}.
\newblock {\em Journal of Physical Chemistry}, 81(25):2340--2361, 1977.

\bibitem{Gil92}
D.~T. Gillespie.
\newblock {A rigorous derivation of the chemical master equation}.
\newblock {\em Physica A: Statistical Mechanics and its Applications},
  188(1-3):404--425, 1992.

\bibitem{Gil01}
D.~T. Gillespie.
\newblock {Approximate accelerated stochastic simulation of chemically reacting
  systems}.
\newblock {\em Journal of Chemical Physics}, 115(4):1716--1733, 2001.

\bibitem{Gil02}
D.~T. Gillespie.
\newblock {The chemical Langevin and Fokker-Planck equations for the reversible
  isomerization reaction}.
\newblock {\em Journal of Physical Chemistry A}, 106(20):5063--5071, 2002.

\bibitem{GuB10}
G.~Guennebaud and B.~Jacob.
\newblock {Eigen v3}, 2010.
\newblock \urlprefix\url{http://eigen.tuxfamily.org/}.

\bibitem{HMT17}
C.~B. Hammouda, A.~Moraes, and R.~Tempone.
\newblock {Multilevel hybrid split-step implicit tau-leap}.
\newblock {\em Numerical Algorithms}, 74(2):527--560, 2017.

\bibitem{HVH19}
X.~Han, M.~Valorani, and N.~N. Habib.
\newblock {Explicit time integration of the stiff chemical Langevin equations
  using computational singular perturbation}.
\newblock {\em Journal of Chemical Physics}, 150:194101, 2019.

\bibitem{HaR02}
E.~L. Haseltine and J.~B. Rawlings.
\newblock {Approximate simulation of coupled fast and slow reactions for
  stochastic chemical kinetics}.
\newblock {\em Journal of Chemical Physics}, 117(15):6959--6969, 2002.

\bibitem{HZC97}
J.~He, H.~Zhang, J.~Chen, and Y.~Yang.
\newblock {Monte Carlo simulation of kinetics and chain length distributions in
  living free-radical polymerization}.
\newblock {\em Macromolecules}, 30(25):8010--8018, 1997.

\bibitem{Hig08}
D.~J. Higham.
\newblock {Modeling and simulating chemical reactions}.
\newblock {\em Siam Review}, 50(2):347--368, 2008.

\bibitem{HAL12}
Y.~Hu, A.~Abdulle, and T.~Li.
\newblock {Boosted hybrid method for solving chemical reaction systems with
  multiple scales in time and population size}.
\newblock {\em Communications in Computational Physics}, 12(4):981--1005, 2012.

\bibitem{Leb94}
V.~I. Lebedev.
\newblock {How to solve stiff systems of differential equations by explicit
  methods}.
\newblock In {\em Numerical methods and applications}, pages 45--80. CRC, Boca
  Raton, FL, 1994.

\bibitem{LeM94}
V.~I. Lebedev and A.~A. Medovikov.
\newblock {Explicit methods of second order for the solution of stiff systems
  of ODEs}.
\newblock {\em Russian Academy of Science}, 1994.

\bibitem{LAE08}
T.~Li, A.~Abdulle, and W.~E.
\newblock {Effectiveness of implicit methods for stiff stochastic differential
  equations}.
\newblock {\em Communications in Computational Physics}, 3(2):295--307, 2008.

\bibitem{Lin72}
B.~Lindberg.
\newblock {IMPEX: a program package for solution of systems of stiff
  differential equations}.
\newblock Technical report, Dept. of Information Processing, Royal Inst. of
  Tech., Stockholm, 1972.

\bibitem{Qua67}
D.~A. McQuarrie.
\newblock {Stochastic Approach to Chemical Kinetics}.
\newblock {\em Journal of Applied Probability}, 4(3):413--478, 1967.

\bibitem{Med98}
A.~A. Medovikov.
\newblock {High order explicit methods for parabolic equations}.
\newblock {\em BIT Numerical Mathematics}, 38(2):372--390, 1998.

\bibitem{RaA03}
C.~V. Rao and A.~P. Arkin.
\newblock {Stochastic chemical kinetics and the quasi-steady-state assumption:
  Application to the Gillespie algorithm}.
\newblock {\em Journal of Chemical Physics}, 118(11):4999--5010, 2003.

\bibitem{RaS07}
M.~Rathinam and H.~{El Samad}.
\newblock {Reversible-equivalent-monomolecular tau: A leaping method for "small
  number and stiff" stochastic chemical systems}.
\newblock {\em Journal of Computational Physics}, 224(2):897--923, 2007.

\bibitem{RPC03}
M.~Rathinam, L.~R. Petzold, Y.~Cao, and D.~T. Gillespie.
\newblock {Stiffness in stochastic chemically reacting systems: The implicit
  tau-leaping method}.
\newblock {\em Journal of Chemical Physics}, 119(24):12784--12794, 2003.

\bibitem{RKV19}
V.~Reshniak, A.~Khaliq, and D.~Voss.
\newblock {Slow-scale split-step tau-leap method for stiff stochastic chemical
  systems}.
\newblock {\em Journal of Computational and Applied Mathematics}, 361:79--96,
  2019.

\bibitem{SSV98}
B.~P. Sommeijer, L.~Shampine, and J.~G. Verwer.
\newblock {RKC: An explicit solver for parabolic PDEs}.
\newblock {\em Journal of Computational and Applied Mathematics},
  88(2):315--326, 1998.

\bibitem{Tal90}
D.~Talay.
\newblock {Second-order discretization schemes of stochastic differential
  systems for the computation of the invariant law}.
\newblock {\em Stochastics and Stochastic Reports}, 29(1):13--36, 1990.

\bibitem{TaT90}
D.~Talay and L.~Tubaro.
\newblock {Expansion of the global error for numerical schemes solving
  Stochastic Differential Equations}.
\newblock {\em Stochastic Analysis and Applications}, 8(4):483--509, 1990.

\bibitem{TiB04}
T.~Tian and K.~Burrage.
\newblock {Binomial leap methods for simulating stochastic chemical kinetics}.
\newblock {\em Journal of Chemical Physics}, 121(21):10356--10364, 2004.

\bibitem{HoS80}
P.~J. {Van der Houwen} and B.~P. Sommeijer.
\newblock {On the internal stability of explicit, $m$-stage Runge--Kutta
  methods for large $m$-values}.
\newblock {\em Zeitschrift f{\"{u}}r Angewandte Mathematik und Mechanik},
  60(10):479--485, 1980.

\bibitem{Ver80}
J.~G. Verwer.
\newblock {An implementation of a class of stabilized explicit methods for the
  time integration of parabolic equations}.
\newblock {\em ACM Transactions on Mathematical Software (TOMS)},
  6(2):188--205, 1980.

\bibitem{VHS90}
J.~G. Verwer, W.~Hundsdorfer, and B.~P. Sommeijer.
\newblock {Convergence properties of the Runge--Kutta--Chebyshev method}.
\newblock {\em Numerische Mathematik}, 57(1):157--178, 1990.

\bibitem{Vil15}
G.~Vilmart.
\newblock {Postprocessed integrators for the high order integration of ergodic
  SDEs}.
\newblock {\em SIAM Journal on Scientific Computing}, 37(1):A201--A220, 2015.

\bibitem{YRS11}
Y.~Yang, M.~Rathinam, and J.~Shen.
\newblock {Integral tau methods for stiff stochastic chemical systems}.
\newblock {\em Journal of Chemical Physics}, 134(4):044129, 2011.

\end{thebibliography}
}{}

\section{Conclusion}\label{sec:conclu}
Based on stabilized methods, second kind Chebyshev polynomials and a postprocessing procedure, we have proposed a new explicit $\tau$-leap method for discrete ergodic stochastic systems with multiple scales (\cref{algo:pseudocode}). Robustness (accuracy of the scheme) and extended stability domains growing quadratically with the stage number have been established (\cref{thm:sktaurock}). Accurate approximation of the invariant measure of ergodic systems has been shown (\cref{thm:post}) thanks to the cheap postprocessing procedure.
Compared to other existing methods, the PSK-$\tau$-ROCK method is shown to be:
\begin{itemize}
	\item faster,
	\item more accurate,
	\item easier to implement.
\end{itemize}
Numerical experiments confirmed the theoretical stability and accuracy results illustrating the efficiency of the PSK-$\tau$-ROCK scheme when compared to other $\tau$-leap schemes for stiff and ergodic discrete stochastic systems: in all considered cases the PSK-$\tau$-ROCK method was the fastest and most accurate scheme.

\section*{Acknowledgments} The authors are partially supported by the Swiss National Science Foundation, under grant No. $200020\_172710$.

%\bibliographystyle{habbrv}
%\bibliography{../../../../../LaTeX/library}

\end{document}